\documentclass[moor]{informs1} 

\usepackage[utf8]{inputenc}
\usepackage[english]{babel}
\usepackage{textcomp}
\usepackage{mathtools}
\usepackage{enumerate}
\usepackage{accents}

\usepackage{natbib}
 \NatBibNumeric
 \def\bibfont{\small}%
 \def\newblock{\ }%
 \bibpunct[, ]{[}{]}{,}{n}{}{,}%

\usepackage[colorlinks=true,breaklinks=true,bookmarks=true,urlcolor=blue,
     citecolor=blue,linkcolor=blue,bookmarksopen=false,draft=false]{hyperref}

\def\EMAIL#1{\href{mailto:#1}{#1}}

\TheoremsNumberedThrough     

\EquationsNumberedThrough    

\newtheorem{axiom}{Axiom}

\begin{document}

\TITLE{\begin{LARGE} Projective functions \end{LARGE}}

\ARTICLEAUTHORS{%
\AUTHOR{Laurence Carassus}
\AFF{ MICS, Centrale-Sup\'{e}lec, Universit\'{e} Paris-Saclay, France
 \EMAIL{laurence.carassus@centralesupelec.fr}}
\AUTHOR{Massinissa Ferhoune}
\AFF{Laboratoire de Math\'{e}matiques de Reims, UMR9008 CNRS and Universit\'{e} de Reims
Champagne-Ardenne, France, \EMAIL{fmassinissa@free.fr}}
} 
\ABSTRACT{
We study the analytical properties of projective functions. We prove that projective functions generalise lower and upper-semianalytic ones while being stable by composition and difference.
 We also show that the class of projective functions is closed under sums,  products, finite suprema and infima, and sections. Assuming the set-theoretical axiom of Projective Determinacy, we prove measurable selection results, stability under integration, and the existence of $\epsilon$-optimal selectors. Finally, we illustrate how these results are important in the context of model uncertainty. }

\KEYWORDS{projective functions; projective sets;  (PD) axiom }
\MSCCLASS{ 03E75, 28A05, 28A25, 28E15, 28B20, 03E60, 54H05, 91G80}

\maketitle


\section{Introduction}
Motivated by model uncertainty in mathematical finance and economics, this paper studies the class of projective functions. 
The definition of projective functions is based on projective sets, which are constructed recursively starting with the Borel sets. Let  $X$ be some Polish space and ${\cal N}$ the Baire space. The analytic sets of $X$ are the projections on $X$ of some Borel sets of $X\times {\cal N}$ and the coanalytic sets are the complements of analytic sets. The class of coanalytic sets is not closed under projection (and complement). So, an analytic set of level 2 is the projection on $X$ of some coanalytic set of $X\times {\cal N}$, while the coanalytic sets of level 2 are the complements of the analytic sets of level 2. But, again, the class of coanalytic sets of level 2 is not closed under projection. 
So, the class of analytic sets of level $n$, $\Sigma_n^1(X)$, and its complement, the class of coanalytic sets of level $n$, $\Pi_n^1(X)$, are constructed recursively, and we denote by  $\Delta_n^1(X)$ their intersection. A set is called projective if it belongs to $\Delta_n^1(X)$ for some $n$, and we denote the class of projective sets by $\mathbf{P}(X)$. \\
We now quickly recall the history of projective sets, for which we refer to the introduction to Moschovakis' manual, see \cite{Moschovakis}. 
Analytic sets were introduced by Suslin in 1917, who realised that the projection of a Borel set is not always a Borel set. Then, Lusin and Sierpinski introduced projective sets in 1925. 
As Moschovakis says \lq\lq{}it was clear from the very beginning that the theory of projective sets was not easy\rq\rq{}. There was two questions of that time which are of interest in this paper: \lq\lq{}Are  $\Sigma_2^1$ sets Lebesgue measurable?\rq\rq{} and \lq\lq{}For a set $A$ in a given class, is it possible to find a function $f$ which graph is included in $A$ and is in the same class that $A$?\rq\rq{} The second question is called the uniformisation problem, and Kondo showed in 1938 that the class of coanalytic sets ($\Pi_1^1$) and the class of analytic sets of level 2 ($\Sigma_2^1$) can be uniformised in the same classes. The next step was made by researchers in descriptive set theory. 
G{\"o}del announced in 1938 that there exists a $\Sigma_2^1$ set, in the model $L$ of constructible sets,  that is not Lebesgue measurable. So, in the Zermelo-Fraenkel set theory with the Axiom of Choice (ZFC), it is not possible to prove that all $\Sigma_2^1$ (and also projective) sets are Lebesgue measurable. Then, Solovay showed in 1970, using forcing methods, that one can consistently assume the axioms
of ZFC with the proposition that all projective sets are Lebesgue measurable. Both results prove that in ZFC, we can neither prove nor disprove the Lebesgue measurability of $\Sigma_2^1$ sets. 
Then, Levy in 1965 showed that we cannot prove in ZFC that every $\Pi_2^1$ set in ${\cal N} \times {\cal N}$ can be uniformised by some projective set. 
The next step was to use strong set-theoretic hypotheses to prove significant results about projective sets. Here we just speak about determinacy. 
Determinacy refers to the existence of a winning strategy for one of the two players of an infinite game. Gale and Stewart showed in 1953 that open (or closed set) are determined. 
Then, Martin showed in 1975 that Borel sets are also determined. 
However, the determinacy of Borel sets is the best possible result provable in ZFC. This is why in 1968, Martin and Moschovakis independently introduced the axiom of Projective Determinacy (PD). The (PD) axiom states that every projective set is determined. Moschovakis showed in 1971 that under this assumption, 
every projective set can be uniformised by a function with a projective graph. 
Mycielski and Swierczkowski showed that under the (PD) axiom, every projective set is universally measurable. 
Then, a surprising connection between the (PD) axiom and large cardinal hypotheses was made. 
Martin and Steel (see \cite{reflargedet_pj}) showed that if sufficiently many Woodin cardinals exist, then all projective sets are determined. Conversely, assuming the (PD) axiom, one can construct inner models (transitive set-theoretic classes that satisfy the axioms of ZF and contain all the ordinals) which exhibit sequences of Woodin cardinals (see \cite{refwood1_pj}). 
The relative consistency of large cardinals with ZFC has been extensively studied via inner model theory and forcing techniques, and no contradictions are known assuming that ZFC itself is consistent (see \cite{refwood1_pj}). 
Quoting Woodin \lq\lq{}Projective Determinacy is the correct axiom for the projective sets; the ZFC axioms are obviously incomplete and, moreover, incomplete in a fundamental way.\rq\rq{}


The notion of projective functions was introduced to state the uniformisation theorem. It is a function whose graph is a projective set. 
We were not aware of this definition when we started to work on the subject, motivated by applications of measure theory to economics and finance.  
So, we introduce an alternative definition of projective function in this paper and the companion paper  \cite{cf2024}, which we call projective measurability and which is more in line with the usual definition of measurability in analysis.  
A function $f:X \to Y$ is projectively measurable if there exists some $n$ such that $f$ is $\Delta_n^1(X)$-measurable. 
To the best of your knowledge, the analytical properties of the projective functions have never been studied before. 
We first remark that Borel, lower, and upper-semianalytic functions are projective. Thus, our notion of measurability generalises the lower and upper-semianalytic one. We establish that sums, differences, products, finite suprema or infima, compositions, and sections of projectively measurable functions remain projectively measurable. 
The closeness results about composition and difference are significant as they do not hold for lower and upper-semianalytic functions. The closeness under composition allows us to show that our definition is equivalent to the definition of projective function in descriptive set theory, i.e., that the graph of $f$ belongs to $\Delta_n^1(X\times Y)$ for some $n$.

Unfortunately, there is another \lq{}\lq{}good\rq\rq{} property that projective sets do not share with analytic sets: the countable union or intersection of projective sets may fail to be a projective set. 
So, the countable infimum or supremum of projective functions may not be projective. 
However, we show that the uncountable infimum or supremum of projective functions on the section of a projective set remains projectively measurable. Then, we provide further results that necessitate the (PD) axiom. We show that the uniformisation theorem in a general Polish space: a projective selection can be made on a projective set.  
We also prove that the integral of a projective function with respect to a projectively measurable stochastic kernel is a projective function and that a projectively measurable $\epsilon$-optimal selector exists for optimising a projective function. Finally, we prove that the multiple-priors conditional support of a projective function is a projective set-valued mapping. These three last properties and the stability by composition are the main results of this paper. 

To show properties on projective functions, we first prove them for $\Delta_n^1$-measurable ones. For that, we use that $\Delta_n^1$ is a sigma-algebra that is closed under Borel preimage. However, $\Delta_n^1$ may not be closed under Borel image, so we will also work with the analytic sets of level $n.$ Indeed, as $\Delta_n^1 \subset \Sigma_n^1$ and $\Sigma_n^1$ is closed by Borel image, the projection on $X$ of some set of $\Sigma_n^1(X \times Y)$ belongs to $\Sigma_n^1(X).$ Then, the inclusion  
$\Sigma_n^1  \subset \Delta_{n+1}^1$ allows us to go back to projective sets at another level.  A lot of results are proved by induction, using the already known results for the analytic sets ($\Sigma_1^1$) and also a parametrisation of the Borel sigma-algebra ($\Delta_1^1$) given in \cite{refproj_pj}. Some of the proofs are completely original and differ from the one for analytic sets  
(for example, the ones of  Lemmata \ref{reciprocal_proj_pj}  and \ref{lemma_before_crit_pj} and Propositions    \ref{prop_crutial_hierarchy_pj} and \ref{univ_cvt_pj}). 

Projective sets originated in measure theory, then became an object of study for logicians. This article, motivated by applications in mathematical finance and economics, proposes a way back to measure theory by studying the analytical properties of projective functions and, in particular, integration and $\epsilon$-optimal selector. We will explain how projective functions and projective sets 
provide a powerful foundation for understanding dynamic decision-making under Knightian model ambiguity. They can capture the complexities of uncertainty without the traps coming from measurable selections.  

 The rest of this article is organised as follows. In Section \ref{subsec_def_proj_pj}, we introduce the projective hierarchy. Then, in Section \ref{subsec_proj_mes_def_pj}, we define projective functions and prove the properties that hold without the (PD) axiom. 
 Section \ref{pd_exp_lab} introduces the (PD) axiom, while Section \ref{subsec_csq_pd_pj} gives some critical consequences of this axiom. In Section \ref{subsec_csq_pd2_pj},  we prove the results holding under the (PD) axiom. Section \ref{secappli} presents applications of projective sets and projective functions when dealing with Knightian uncertainty in mathematical finance and economics. 
Finally, the appendix contains the missing proofs.

\section{Projective sets}
\label{subsec_def_proj_pj}
We give the definition of projective sets (see Definition \ref{def_proj_set_pj}). Then, Proposition \ref{base_hierarchy_pj} shows that these sets share a lot of properties with the analytic sets. Projective sets are also closed under complement (in opposition to the class of analytic sets). But still, as for analytic sets, the class of projective sets is not a sigma-algebra: a countable union (or intersection) of projective sets may not be a projective set. This is why the class $\Delta_{n}^1(X)$ is important: it is a sigma-algebra. 
\subsection{Definition of projective sets}
Projective sets can be seen as a generalisation of analytic sets. 
Let $X$ be a Polish space. An analytic set of $X$ is the projection on $X$ of the Borel subsets of $X\times\mathbb{N}^{\mathbb{N}}$, see \citep[Proposition 7.39, p165]{ref1_pj}. We will use this definition of analytic sets throughout the paper. Other equivalent definitions can be found in  \citep[Proposition 7.41, p166]{ref1_pj}. We denote by $\Sigma_1^1(X)$ the class of analytic sets of $X$. 
Then, the projection on $X$ of any set of $\Sigma_1^1(X\times \mathbb{N}^{\mathbb{N}})$ is in $\Sigma_1^1(X)$, see \citep[Proposition 7.39, p165]{ref1_pj}. However, $\Sigma_1^1(X)$ is not closed under complement, see \citep[Proposition B.6, p292]{ref1_pj}. The complement of an analytic set is called a coanalytic set. We denote by $\Pi_1^1(X)$, the class of coanalytic sets of $X$ and by $\Delta_1^1(X)$, the intersection of $\Sigma_1^1(X)$ and $\Pi_1^1(X)$. 
Then, $\Delta_1^1(X)$ is the class of Borel sets, denoted by ${\cal B}(X)$ (see \citep[Theorem 14.11, p88]{refproj_pj}). 
Now, the projection on $X$ of any set of $\Pi_1^1(X\times \mathbb{N}^{\mathbb{N}})$ has no reason to be coanalytic or analytic, and we denote by $\Sigma_2^1(X)$ the class of these sets. The elements of $\Sigma_2^1(X)$ can be seen as analytic sets of level 2. 
Repeating inductively this scheme and taking the (increasing) union of the $\Delta_n^1(X)$ gives the class $\textbf{P}(X)$ of projective sets. 
The reader can find an elaborate construction of them in \citep[p313]{refproj_pj}. 
\begin{definition}
Let $X$ be a Polish space. Let $\mathcal{N}:= \mathbb{N}^\mathbb{N}$ be the Baire space. First, $\Sigma_1^1(X)$ is the class of analytic sets of $X,$ while $\Pi_1^1(X)$ is the class of coanalytic sets of $X$:
\begin{eqnarray}
\Sigma_{1}^1(X) := \{\textup{proj}_X(C),\;  C\in {\cal B}(X\times \mathcal{N})\}\label{Sigma_pj} \quad 
\Pi_{1}^1(X) := \{X\setminus C,\; C\in\Sigma_{1}^1(X)\}. 
\label{level1}
\end{eqnarray}
 For each $n\geq 2$, we define recursively the classes $\Sigma_n^1(X)$ and $\Pi_n^1(X)$ of $X$ as follows:
\begin{eqnarray}
\Sigma_{n}^1(X)&:=& \{\textup{proj}_X(C),\;  C\in\Pi_{n-1}^1(X\times \mathcal{N})\}\label{Sigma_pj}\\
\Pi_{n}^1(X)&:=& \{X\setminus C,\; C\in\Sigma_{n}^1(X)\}, \label{Pi_pj}
\end{eqnarray}
Then, for all $n\geq 1$, we set 
\begin{eqnarray}
\Delta_{n}^1(X)&:=& \Sigma_n^1(X)\cap \Pi_n^1(X) \label{Delta_pj}.
\end{eqnarray}
The class $\mathbf{P}(X)$ of projective sets of $X$ is defined by $$\mathbf{P}(X) := \bigcup_{n\geq 1} \Delta_n^1(X).$$ 
\label{def_proj_set_pj}
\end{definition}
\subsection{Properties of projective sets}
The following proposition justifies the term ``hierarchy" and gives crucial properties of $\Sigma_{n}^1(X)$, $\Pi_n^1(X)$, $\Delta_n^1(X)$ and $\textbf{P}(X).$ 
\begin{proposition}
Let $n\geq 1,$ $X$ and $Y$ be Polish spaces. Let $f : X \to Y$ be Borel measurable.\\
(i) The class $\Sigma_n^1(X)$ is closed under countable intersections and unions. Moreover,  we have that $f^{-1}(B)\in \Sigma_n^1(X)$ for all $B\in\Sigma_n^1(Y)$ and  $f(A)\in\Sigma_n^1(Y)$ for all $A\in\Sigma_n^1(X)$.\\ 
(ii) The class $\Pi_n^1(X)$ is closed under countable intersections and unions. Moreover, we have that $f^{-1}(B)\in \Pi_n^1(X)$ for all $B\in\Pi_n^1(Y)$.\\ 
(iii) The class $\Delta_n^1(X)$ is a sigma-algebra on $X,$ $\Delta_1^1(X)=\mathcal{B}(X)$ and $f^{-1}(B)\in \Delta_n^1(X)$ for all $B\in\Delta_n^1(Y)$.\\ 
(iv) The sequences $(\Sigma_n^1(X))_{n\geq 1}$, $(\Pi_n^1(X))_{n\geq 1}$ and $(\Delta_n^1(X))_{n\geq 1}$ are nondecreasing.\\ 
(v) We have that 
$$\Sigma_n^1(X) \times \Sigma_n^1(Y) \subset \Sigma_n^1(X\times Y), \, \Pi_n^1(X) \times \Pi_n^1(Y) \subset \Pi_n^1(X\times Y), \, \Delta_n^1(X) \times \Delta_n^1(Y) \subset \Delta_n^1(X\times Y)$$ 
\begin{eqnarray}
\Sigma_n^1(X) \cup \Pi_n^1(X) \subset \Delta_{n+1}^1(X)\;\; \mbox{and}\;\; \textbf{P}(X)=\bigcup_{n\geq 1} \Delta_n^1(X)= \bigcup_{n\geq 1} \Pi_n^1(X)=\bigcup_{n\geq 1} \Sigma_n^1(X) . \label{eq(v)_pj}
\end{eqnarray}
(vi)  The class $\textbf{P}(X)$ is closed under finite unions and intersections. It is also closed under complements. Moreover, we have that $f^{-1}(B)\in \textbf{P}(X)$ for all $B\in\textbf{P}(Y)$ and $f(A) \in \textbf{P}(Y)$ for all $A\in\textbf{P}(X)$. Finally, 
\begin{eqnarray}
\mathcal{B}(X)\subset \textbf{P}(X),\;\; \Sigma_1^1(X)\cup \Pi_1^1(X)\subset \textbf{P}(X),\;\; \textbf{P}(X) \times \textbf{P}(Y) \subset \textbf{P}(X\times Y).
\label{eq(vi)_pj}
\end{eqnarray}
\label{base_hierarchy_pj}
\end{proposition}
\proof{Proof.}
See Appendix \ref{hiera_section_proof} .
\Halmos \\ \endproof

\begin{remark}
Proposition \ref{base_hierarchy_pj} is not \citep[Proposition 37.1, p314]{refproj_pj}, as the classes of sets $\Sigma_n^1$ and $\Pi_n^1$ in 
\cite{refproj_pj} differ from ours. Indeed, they are defined by
\begin{eqnarray*}
\Sigma_n^1 := \bigcup_{X \textup{Polish space}} \Sigma_{n}^1(X) \quad \mbox{and} \quad \Pi_n^1 := \bigcup_{X \textup{Polish space}} \Pi_{n}^1(X).
\end{eqnarray*}
\end{remark}

\begin{remark}
The classes $\Pi_n^1(X)$ and $\Delta_n^1(X)$ may not be closed under Borel image. Indeed, (ii) and (iii) follow from (i) as $f^{-1}(Y \setminus C) = X \setminus f^{-1}(C)$ but $f(X\setminus C)$ and $Y\setminus f(C)$ have no reason to be included in each other. 
\end{remark}
\begin{claim}
 A countable union of projective sets may not be a projective set; thus, the class of projective sets may fail to be a sigma-algebra. 
 \end{claim}
 \proof{Proof.}
Let $X$ be an uncountable Polish space and $(A_n)_{n\geq 2}$ be such that $A_n\in\Sigma_n^1(X)$ and $A_n\notin\Sigma_{n-1}^1(X)$ for all $n\geq 2$. Such a sequence exists (see \citep[Theorem 37.7, p316]{refproj_pj} and \eqref{Delta_pj}). \\
Let $B_n := \{(n,x),\; x\in A_n\}$ for all $n\geq 2$. We prove below that $B_n \in \textbf{P}(\mathbb{N}\times X)$ for all $n\geq 2$, but that $B:= \bigcup_{n\geq 2} B_n \notin \textbf{P}(\mathbb{N}\times X)$. Thus, $\textbf{P}(\mathbb{N}\times X)$ is not a sigma-algebra of $\mathbb{N}\times X$.\\ Note first that for all $n\geq 2$, $B_n = \ell_n(A_n)$ where $\ell_n(x):=(n,x)$ for all $x\in X$. 
As $A_n\in\Sigma_n^1(X)$ and $\ell_n$ is Borel measurable, we get that $B_n \in \Sigma_n^1(\mathbb{N}\times X),$ see (i) in Proposition \ref{base_hierarchy_pj}. So, \eqref{eq(v)_pj} in Proposition \ref{base_hierarchy_pj} shows that for all $n\geq 2$
$$B_n \in \Sigma_n^1(\mathbb{N}\times X) \subset \Delta_{n+1}^1(\mathbb{N}\times X),$$  
and $B_n\in   \textbf{P}(\mathbb{N}\times X).$ 
Assume now that $B \in \textbf{P}(\mathbb{N}\times X)$. Then, there exists $p\geq 1$ such that $B\in \Delta_p^1(\mathbb{N}\times X) \subset \Sigma_p^1(\mathbb{N}\times X)$ (see \eqref{Delta_pj}). Remark that 
\begin{eqnarray*}
\ell_{p+1}^{-1}(B) = \{x\in X,\; (p+1,x)\in B\} &=& \bigcup_{n\geq 2}\{x\in X,\; (p+1,x)\in B_n\} = \{x\in X,\; (p+1,x)\in B_{p+1}\} = A_{p+1}.
\end{eqnarray*}
As $B\in\Sigma_p^1(\mathbb{N}\times X)$ and $\ell_{p+1}$ is Borel measurable, using again (i) in Proposition \ref{base_hierarchy_pj}, $A_{p+1}\in \Sigma_p^1(X)$, a contradiction. So, $B\notin \textbf{P}(\mathbb{N}\times X)$.
\Halmos \\ \endproof

\section{Projective functions}
 \label{subsec_proj_mes_def_pj}
\subsection{Definition of projective functions} 
 As mention in the introduction, we were not aware of the definition of projective functions as functions whose graph is a projective set used, for example, in \cite{woodincs}.   
So, we introduce an alternative definition of projective function, which we call projective measurability, and is more in line with the usual definition of measurability in analysis.  
We will show that our definition coincides with the one in descriptive set theory. 
The key concept of projective functions will allow us to generalise properties that hold for upper and lower semi-analytic functions.\\
\begin{definition}
Let $X$ and $Y$ be Polish spaces and $D\subset X$. A function $f : D \to Y$ is projective (or projectively measurable) if $D\in\textbf{P}(X)$ and if there exists some $n\geq 1$ such that $f$ is $\Delta_n^1(X)$-measurable in the sense that for all $B\in\mathcal{B}(Y)$, $f^{-1}(B)=\{x\in D, f(x)\in B\}\in\Delta_n^1(X)$. \\
\label{proj_fct_def_pj}
\end{definition}
\begin{remark}
(i) It is important in Definition \ref{proj_fct_def_pj} that the $n$ is the same for all  $B\in\mathcal{B}(Y)$. \\
(ii) We found two references to projective functions with a similar definition in an exercise of the textbook of Kechris, see \cite{refproj_pj}. However, this notion is not exploited, and no results are proved. \\
(iii) As $\Delta_n^1(X)$ is a sigma-algebra (see (iii) in Proposition \ref{base_hierarchy_pj}), to prove that $f : X \to \mathbb{R}\cup \{-\infty,+\infty\}$ is $\Delta_n^1(X)$-measurable (resp. projective), it is enough to prove that $\{f<c\}$ or $\{f\leq c\}$ belong to $\Delta_n^1(X)$ for all $c\in\mathbb{R}$ (resp. for some given $n\geq1$). This will be used in the rest of the paper without further mention.\\
(iv) To show properties on projective functions, we will first prove them for $\Delta_n^1(X)$-measurable ones. For that, we either use that $\Delta_n^1(X)$ is a sigma-algebra,  stable by Borel preimage, or if we need the direct image, that $\Sigma_n^1(X)$ is stable under Borel image. Then, \eqref{eq(v)_pj} shows that 
$\Sigma_n^1(X)  \subset \Delta_{n+1}^1(X)$ and we can go back to some $\Delta_{p}^1(X)$-measurability at another level of the hierarchy and also to projective measurability. \\
\label{remarkxx}
\end{remark}
As a direct consequence of Definition \ref{proj_fct_def_pj} and Proposition \ref{base_hierarchy_pj}, we find that projective functions generalise Borel and upper and lower-semianalytic ones. 
\begin{corollary}
Let $X$ and $Y$ be Polish spaces and let $f: X\to Y$.\\ (i) If $f$ is $\mathcal{B}(X)$-measurable, then $f$ is a projective function.\\ (ii) If $Y= \mathbb{R}\cup\{-\infty, +\infty\}$ and $f$ is a lower-semianalytic  or upper-semianalytic function (see \citep[Definition 7.21, p177]{ref1_pj}), then $f$ is a projective function.
\label{borel_proj_pj}
\end{corollary}
\proof{Proof.}
(i) The first assertion follows from $\mathcal{B}(X) =\Delta_1^1(X)$, see  (iii) in Proposition \ref{base_hierarchy_pj}. \\
(ii) Let $c\in\mathbb{R}$. If $f$ is lower (resp. upper)-semianalytic, $\{f<c\}$ (resp. $\{f>c\}$) belongs to $\Sigma_1^1(X)$ and we conclude as $\Sigma_1^1(X)\subset \Delta_2^1(X)$ (see \eqref{eq(v)_pj}). 
 \Halmos \\ \endproof

\subsection{Property of general projective functions} 
We now prove properties of general projective functions, which are true without the (PD) axiom. We show that the composition of projective functions remains projective. Note that the universally measurable functions also have this property (see \citep[Proposition 7.44, p172]{ref1_pj}) but not the lower and upper-semianalytic ones. 
We also show that the sections of projective functions are projective. 
The first two lemmata are simple and will be valuable tools for several proofs. The first one shows that a vector of projective functions is projective.
\begin{lemma}
Let $X$, $Y$ and $Z$ be Polish spaces and let $f: X\to Y$ and $g : X \to Z$. Let $h : X \to Y \times Z$ be defined by $h(x) := (f(x),g(x))$ for all $x\in X$. \\
(i) If $f$ and $g$ are $\Delta_p^1(X)$-measurable for some $p\geq 1$, then $h$ is $\Delta_p^1(X)$-measurable.\\
(ii) If $f$ and $g$ are projective functions, then $h$ is also projective.
\label{couple_fct_mes}
\end{lemma}
\proof{Proof.}
(i) Let $A\in\mathcal{B}(Y)$ and $B\in \mathcal{B}(Z)$. We have that $h^{-1}(A\times B) = f^{-1}(A) \cap g^{-1}(B)$ and that $f^{-1}(A)$ and $g^{-1}(B)$ belong to $\Delta_{p}^1(X)$. Thus, Proposition \ref{base_hierarchy_pj} (iii) ensures that $h^{-1}(A\times B) \in \Delta_p^1(X)$ and $h$ is $\Delta_p^1(X)$-measurable. \\
(ii) If $f$ and $g$ are  projective functions, there exist $n, p \geq 1$ such that $f$ is $\Delta_n^1(X)$-measurable and $g$ is $\Delta_p^1(X)$-measurable. We may assume that $n\leq p$. Then, Proposition \ref{base_hierarchy_pj} (iv) shows that $f$ is $\Delta_p^1(X)$-measurable, and (ii) follows from (i).
\Halmos\\ \endproof

\begin{lemma}
Let $X$, $Y$ and $Z$ be Polish spaces and let $f : X \to Z$. Let $h : X  \times Y \to  Z$ be defined by $h(x,y) := f(x)$ for all $(x,y)\in X\times Y$. \\
(i) If $f$ is $\Delta_p^1(X)$-measurable for some $p\geq 1$, then $h$ is $\Delta_p^1(X \times Y)$-measurable.\\
(ii) If $f$ is a projective function, then $h$ is also projective.
\label{etendre_fct_mes}
\end{lemma}
\proof{Proof.}
(i) Let $B\in \mathcal{B}(Z)$. We have that 
$$h^{-1}(B) = \{(x,y)\in X\times Y, \, f(x) \in B\}=f^{-1}(B) \times Y\in \Delta_{p}^1(X)\times \Delta_{p}^1(Y) \subset \Delta_{p}^1(X \times Y),$$ 
using Proposition \ref{base_hierarchy_pj} (iii)-(v) as $Y \in {\cal B}(Y)=\Delta_{1}^1(Y)$. So,   $h$ is $\Delta_p^1(X \times Y)$-measurable. \\
(ii) If $f$ is a projective function, there exists $p\geq 1$ such that $f$ is $\Delta_p^1(X)$-measurable. So, from (i) $h$ is $\Delta_p^1(X \times Y)$-measurable and thus projective. 
\Halmos\\ \endproof

We now extend (iii) in Proposition \ref{base_hierarchy_pj} to $\Delta_p^1(X)$-measurable functions. This result is a first step to proving the closeness under composition and also that our definition coincides with the one in descriptive set theory.
\begin{lemma}
Let $X$ and $Y$ be Polish spaces. Let $p\geq 1$ and assume that $f : X \to Y$ is $\Delta_p^1(X)$-measurable. Then, for all $n\geq 1$ and for all $B \in \Delta_{n}^1(Y)$, $f^{-1}(B)\in  \Delta_{p+n}^1(X)$.
\label{reciprocal_proj_pj}
\end{lemma}
\proof{Proof.}
First, we show by induction on $n$ that for all $n\geq 1$, $p\geq 1$, for all Polish spaces $X$ and $Y$, for all $\Delta_p^1(X)$-measurable $f : X\to Y$  and for all $B\in \Sigma_{n}^1(Y)$, we have that $f^{-1}(B) \in \Sigma_{n+p-1}^1(X)$. We begin with the heredity step. Assume that the induction hypothesis holds for $n\geq 1$. Let $B\in\Sigma_{n+1}^1(Y)$. There exists some $C\in \Pi_{n}^1(Y\times \mathcal{N})$ such that $B= \textup{proj}_Y(C)$. Let $p\geq 1$ and  $f : X \to Y$ be $\Delta_p^1(X)$-measurable.  For all $(x,u)\in X\times\mathcal{N}$, let $\Psi(x,u):=(f(x),u)\in Y\times \mathcal{N}.$ Then, 
\begin{eqnarray}
f^{-1}(B) = \{x\in X,\; f(x) \in B  \} &=&\{x\in X,\; \exists u\in\mathcal{N},\; (f(x),u) \in C \}\nonumber\\ &= & \{x\in X,\; \exists u\in\mathcal{N},\; (x,u) \in \Psi^{-1}(C)\}=\textup{proj}_X(\Psi^{-1}(C)). \label{temp_eq_set_recipq_pj}  
\end{eqnarray}
Note that $(x,u) \mapsto f(x)$ and $(x,u) \mapsto u$ are $\Delta_p^1(X\times \mathcal{N})$-measurable using (iii) to (v) in Proposition \ref{base_hierarchy_pj}. So, Lemma \ref{couple_fct_mes} shows that $\Psi$ is $\Delta_p^1(X\times \mathcal{N})$-measurable. As $C\in \Pi_n^1(Y\times \mathcal{N})$, \eqref{Pi_pj} implies that $(Y\times \mathcal{N})\setminus C \in\Sigma_{n}^1(X\times \mathcal{N})$ and using the induction hypothesis for $X\times \mathcal{N},$ $\Psi$ and $(Y\times \mathcal{N})\setminus C$, we get that 
$$(X\times \mathcal{N})\setminus \Psi^{-1}(C)=\Psi^{-1}((Y\times \mathcal{N})\setminus C) \in \Sigma_{n+p-1}^1(X\times \mathcal{N}).$$ 
So,  $\Psi^{-1}(C)\in \Pi_{n+p-1}^1(X\times \mathcal{N})$. Finally,   \eqref{temp_eq_set_recipq_pj}   and \eqref{Sigma_pj} show that $f^{-1}(B) \in \Sigma_{n+p}^1(X)$.\\
We now turn to the initialisation step. Let $B\in\Sigma_1^1(Y)$. 
There exists $C \in \mathcal{B}(Y\times \mathcal{N})$ such that $B= \textup{proj}_Y(C)$, see \eqref{level1}. Let $p\geq 1,$  $f : X \to Y$ be $\Delta_p^1(X)$-measurable and $\Psi$ be defined as above. 
Then, again $\Psi$ is $\Delta_p^1(X\times \mathcal{N})$-measurable and  \eqref{Delta_pj} shows that
$$\Psi^{-1}(C) \in \Delta_p^1(X\times \mathcal{N}) \subset \Sigma_p^1(X\times \mathcal{N}).$$ 
Recalling \eqref{temp_eq_set_recipq_pj}, $f^{-1}(B) =\textup{proj}_X(\Psi^{-1}(C))$.  
So, as $\textup{proj}_X$ is Borel measurable, (i) in Proposition \ref{base_hierarchy_pj} shows that $f^{-1}(B) \in \Sigma_{p}^1(X)$. This concludes the induction.\\ Let $n\geq 1$ and $B \in \Delta_n^1(Y) \subset \Sigma_n^1(Y)$. Let $p\geq 1$ and $f : X \to Y$ be $\Delta_p^1(X)$-measurable. Then, $f^{-1}(B) \in \Sigma_{n+p-1}^1(X) \subset \Delta_{n+p}^1(X)$ (see \eqref{eq(v)_pj}),  which is the desired result.
\Halmos \\\endproof

We now prove that our definition of projective functions is the same as the one used in descriptive set theory, see for example \cite{woodincs}. 
\begin{proposition}
Let $X$ and $Y$ be Polish spaces and $D\subset X$. Let $f : D \to Y.$ \\
(i) Assume that $D \in \Delta_{n}^1(X)$  and $\textup{Graph}(f) \in \Delta_{n}^1(X \times Y).$ Then, $f$ is $\Delta_{n+1}^1(X)$-measurable. \\
(ii) Assume that $D \in \Delta_{n}^1(X)$  and  $f$ is $\Delta_{n}^1(X)$-measurable.  Then, $\textup{Graph}(f) \in \Delta_{n+1}^1(X \times Y).$\\
(iii) Let $D\in\textbf{P}(X)$. The function $f$ is projective if and only if $\textup{Graph}(f) \in \textbf{P}(X \times Y).$
\label{pareil}
\end{proposition}
\proof{Proof.}
(i) Let $B \in {\cal B}(Y).$ Then, 
\begin{eqnarray*}
f^{-1}(B) = \{x\in D,\; f(x) \in B  \} = 
\{x\in X,\; \exists y \in Y,\; y=f(x), \, x \in D, \,  y \in B\}=\textup{proj}_X \big(\textup{Graph}(f) \cap (D \times B)\big).
\end{eqnarray*}
Now,  (iii)  to (v) in Proposition \ref{base_hierarchy_pj} show that $D \times B$ and  $\textup{Graph}(f) \cap (D \times B)$ belong to  $\Delta_{n}^1(X \times Y)$. 
But, $\Delta_{n}^1(X \times Y) \subset \Sigma_{n}^1(X \times Y)$, and as $\textup{proj}_X$ is Borel measurable, (i) in Proposition \ref{base_hierarchy_pj} and \eqref{eq(v)_pj}   show that 
$$f^{-1}(B) \in \Sigma_{n}^1(X)  \subset \Delta_{n+1}^1(X).$$
(ii) The function $ X \times Y \ni (x,y) \mapsto\Phi(x,y):=(f(x),y)\in Y \times Y$ is $\Delta_{n}^1(X \times Y)$-measurable (see Lemma \ref{couple_fct_mes}). Now, 
$$ \textup{Graph}(f) = \{ (x,f(x)), \, x\in D\}=\Phi^{-1}(\{(y,y), \, y \in Y\})\cap (D \times Y).$$
As $\{(y,y), \, y \in Y\}$ is closed, it belongs to ${\cal B}(X \times Y)= \Delta_{1}^1(X \times Y)$. So, Lemma \ref{reciprocal_proj_pj} shows that 
$\Phi^{-1}(\{(y,y), \, y \in Y\}) \in \Delta_{n+1}^1(X \times Y)$. Now (iii)  to (v) in Proposition \ref{base_hierarchy_pj} show that $D \times Y$ and also 
$\textup{Graph}(f)$ belong to $\Delta_{n+1}^1(X \times Y)$.\\
(iii) Let $D\in\textbf{P}(X)$.  There exists some $p\geq 1$ such that $D \in \Delta_p^1(X)$ and for every $ q \geq p$ $D \in \Delta_q^1(X)$ (see (iv) in Proposition \ref{base_hierarchy_pj}).\\
Assume that $f$ is projective. Then,  there exists some $q\geq 1$ such that  $f$ is $\Delta_q^1(X)$-measurable. We may assume that $p\leq q$ and (ii) shows that 
$$\textup{Graph}(f) \in \Delta_{q+1}^1(X \times Y) \subset \textbf{P}(X \times Y).$$
Now if $\textup{Graph}(f) \in \textbf{P}(X \times Y),$  there exists some  $q\geq 1$ such that 
$\textup{Graph}(f) \in \Delta_q^1(X\times Y).$ Again, we may assume that $p\leq q$ and (i) shows that $f$ is $\Delta_{q+1}^1(X)$-measurable and thus projective. 
\Halmos \\\endproof

We can now show that the class of projective functions is closed under composition. 
\begin{proposition}
Let $X, Y$, and $Z$ be Polish spaces. Let  $g: D \to Y$ and $f : E \to Z$ where $D\subset X$ and $g(D)\subset E \subset Y$.\\ (i) Assume that $f$ is $\Delta_p^1(Y)$-measurable and that $g$ is $\Delta_q^1(X)$-measurable for some $p, q\geq 1$. Then, $f\circ g$ is $\Delta_{p+q}^1(X)$-measurable.\\ 
(ii) Assume that $f$ and $g$ are projective functions. Then, $f\circ g$ is also  a projective function.
\label{comp_proj_pj}
\end{proposition}
\proof{Proof.}
We show (i). Let $B\in\mathcal{B}(Z)$. We have that $(f\circ g)^{-1}(B) = g^{-1}(f^{-1}(B))$. As $f$ is $\Delta_p^1(Y)$-measurable, $f^{-1}(B)\in \Delta_p^1(Y)$ and Lemma \ref{reciprocal_proj_pj} shows that $g^{-1}(f^{-1}(B)) \in \Delta_{p+q}^1(X)$: (i) is proved. If $f$ and $g$ are projective functions,  there exist some $p\geq 1$ and $q\geq 1$ such that $f$ is $\Delta_p^1(X)$-measurable and $g$ is $\Delta_q^1(X)$-measurable and (i) shows that $f\circ g$ is $\Delta_{p+q}^1(X)$-measurable and thus projectively measurable : (ii) is proved. 
\Halmos \\\endproof

As a corollary, we obtain that the sections of projective functions are projective.
\begin{corollary}
Let $X, Y$, and $Z$ be Polish spaces. Let $h : X\times Y \to Z$.\\ 
(i)  Assume that $h$ is $\Delta_{p}^1(X\times Y)$-measurable for some $p\geq 1$. Then $h(x,\cdot) : y \mapsto h(x,y)$ is $\Delta_{p+1}^1(Y)$-measurable for all $x\in X$ and $h(\cdot,y) : x \mapsto h(x,y)$ is $\Delta_{p+1}^1(X)$-measurable for all $y\in Y$.\\ 
(ii) Assume that $h$ is a projective function. Then $h(x,\cdot) : y \mapsto h(x,y)$ is projectively measurable for all $x\in X$ and $h(\cdot,y) : x \mapsto h(x,y)$ is projectively measurable for all $y\in Y$.
\label{comp_sec_pj}
\end{corollary}
\proof{Proof.}
We show (i).  Fix $x\in X$ and let $\sigma : y \mapsto (x,y)$. Then, $\sigma$ is Borel measurable and also $\Delta_1^1(Y)$-measurable by Proposition \ref{base_hierarchy_pj} (iii). As $h(x,\cdot) = h \circ \sigma$ and $h$ is $\Delta_p^1(X\times Y)$-measurable, we deduce from Proposition \ref{comp_proj_pj} (i) that $h(x,\cdot)$ is $\Delta_{p+1}^1(Y)$-measurable. The proof that $h(\cdot,y)$ is $\Delta_{p+1}^1(X)$-measurable for all $y\in Y$ is similar and omitted. \\
We show (ii). If $h$ is a projective function, there exists some $p\geq 1$ such that $h$ is $\Delta_p^1(X\times Y)$ and (i) shows that $h(x,\cdot)$ and $h(\cdot,y)$ are respectively $\Delta_{p+1}^1(Y)$ and $\Delta_{p+1}^1(X)$-measurable for all $x\in X$ and $y\in Y$, and thus projectively measurable.
\Halmos \\\endproof

\subsection{Properties of real or vector-valued projective functions}
We now turn to results that hold for real or vector-valued functions. 
The first lemma shows that the set of finite vector-valued projective functions is closed under usual operations.
\begin{lemma}
Let $X$ be a Polish space and let $f: X\to \mathbb{R}^n$ and $g : X \to \mathbb{R}^n$ for some $n\geq 1$. \\
(i) If $f$ and $g$ are $\Delta_p^1(X)$-measurable for some $p\geq 1$, then $f+g$, $\langle f, g \rangle$ and $-f$ are $\Delta_p^1(X)$-measurable. \\
(ii) If $f$ and $g$ are projective functions, then $f+g$, $\langle f, g \rangle$ and $-f$  are also projective functions.
\label{real_vec_fct_mes}
\end{lemma}
\proof{Proof.}
(i) Let $B_1 \in \mathcal{B}(\mathbb{R}^d)$, $B_2 \in \mathcal{B}(\mathbb{R}),$  $B_3\in\mathcal{B}(\mathbb{R}^d)$ and $h(x):=(f(x),g(x))$. Lemma \ref{couple_fct_mes} shows that $h$ is $\Delta_p^1(X)$-measurable. Moreover, we have that
\begin{eqnarray*}
(f+g)^{-1}(B_1) &=& \{x\in X,\; f(x)+g(x) \in B_1\} = \{x\in X,\; (f(x),g(x)) \in \sigma_1^{-1}(B_1)\}=h^{-1}\big(\sigma_1^{-1}(B_1)\big)\\
\langle f, g \rangle^{-1}(B_2) &=& \{x\in X,\; \langle f(x), g(x) \rangle \in B_2\} = \{x\in X,\; (f(x),g(x)) \in \sigma_2^{-1}(B_2)\}=h^{-1}\big(\sigma_2^{-1}(B_2)\big)\\
(-f)^{-1}(B_3) &=& \{x\in X,\; -f(x)\in B_3\} = \{x\in X,\; f(x) \in \sigma_3^{-1}(B_3)\}=f^{-1}\big(\sigma_3^{-1}(B_3)\big),
\end{eqnarray*}
where $\sigma_1(x,y) := x+y$, $\sigma_2(x,y) := \langle x, y \rangle$, $\sigma_3(x) := -x$ for all $(x,y)\in\mathbb{R}^d \times \mathbb{R}^d$. As for all $i\in \{1,2,3\}$, 
 $\sigma_i$ is Borel measurable, $\sigma_i^{-1}(B_i)$ is a borelian, and (i) follows. \\
(ii) This is the same proof as in Lemma \ref{couple_fct_mes} (ii). 
\Halmos\\ \endproof

The next proposition extends the results of Lemma \ref{real_vec_fct_mes} to real but possibly infinite-valued functions under conventions \eqref{cvt_inf_pj} and \eqref{cvt_new} below. 
\begin{eqnarray}
\label{cvt_inf_pj}
-\infty+\infty=+\infty-\infty=-\infty.
\end{eqnarray}
This is the usual convention for the maximisation problems considered in the Application Section \ref{secappli}. 
Otherwise, we adopt the usual arithmetic rules in calculations involving $(+\infty)$ and $(-\infty)$ described in \citep[Section 7.44, p26-27]{ref1_pj} :
\begin{eqnarray}
+(-\infty) = -\infty, \quad -(+\infty) = -\infty \quad 0 \times (\pm \infty) = (\pm \infty) \times 0 = 0. 
\label{cvt_new}
\end{eqnarray}

\begin{proposition} 
Let $X$ be a Polish space and for all $n\geq 0$, let $f, f_n, g : X \to \mathbb{R}\cup \{-\infty,+\infty\}$.\\ 
(i) Let $p \geq 1$. Assume that $f$, $f_n$ and $g$ are $\Delta_p^1(X)$-measurable for all $n\geq 0$. Then, $fg,$ 
$f+g$, 
$-f$, $\min(f,g)$, $\max(f,g)$, $\inf_{n\geq 0} f_n$, $\sup_{n\geq 0} f_n$ are $\Delta_p^1(X)$-measurable. If $f> 0$, $f^a$ is also $\Delta_p^1(X)$-measurable for all $a\neq0$.\\ 
(ii) If $f$ and $g$ are  projective functions, then $fg$, $f+g$, $-f$, $\min(f,g)$ and  $\max(f,g)$  are also projective functions. If $f> 0$, $f^a$ is a projective function for all $a\neq0$.
\label{lemma_real_proj_pj}
\end{proposition}
\proof{Proof.}
See Appendix \ref{hiera_section_proof}.
\Halmos \\ \endproof

Remark that we don't have that if $f_n$ is a projective function for all $n\geq 0$, then $\inf_{n\geq 0} f_n$ or  $\sup_{n\geq 0} f_n$ is projectively measurable. Indeed, each $f_n$ is $\Delta_{p_n}^1(X)$-measurable but we may not find some $p$ such that all the $f_n$ are $\Delta_{p}^1(X)$-measurable. 
However, we show that the uncountable supremum and infimum of projective functions remain projectively measurable if the optimum is taken on the section of a projective set. This proposition can be seen as an extension of \citep[Proposition 7.47, p179]{ref1_pj}, which holds for lower or upper-semianalytic functions. It implies that if $(x,n) \mapsto f(x,n)$ is projective that $\inf_{n\geq 0} f_n$ or  $\sup_{n\geq 0} f_n$ are indeed projective. 
\begin{proposition}
Let $X$ and $Y$ be Polish spaces. Let $D \in \textbf{P}(X\times Y)$ and $f : X \times Y \to \mathbb{R}\cup\{-\infty,+\infty\}$ be a projective function. Let $D_x := \{y\in Y,\; (x,y)\in D\}$ for all $x\in X$. Then, the functions $f_*, f^* : \textup{proj}_X(D) \to \mathbb{R}\cup\{-\infty,+\infty\}$ defined by
\begin{eqnarray*}
f_*(x) := \inf_{y\in D_x} f(x,y) \;\; \mbox{and} \;\; f^*(x) := \sup_{y\in D_x} f(x,y)
\end{eqnarray*}
are projective functions.
\label{prop_infsup_mes_proj_pj}
\end{proposition}

\proof{Proof.}
There exist some $p\geq 1$ and $q\geq 1$ such that $D \in \Delta_p^1(X\times Y)$ and $f$ is $\Delta_q^1(X\times Y)$-measurable. We may assume that $p\leq q$ and $D \in \Delta_q^1(X\times Y)$ (see (iv) in Proposition \ref{base_hierarchy_pj}). 
We show that $f_*$ is a projective function. Let $c\in \mathbb{R}$. We have that 
\begin{eqnarray*}
f_*^{-1}([-\infty,c)) &=& \{x\in X,\; \exists y\in Y,\; (x,y)\in D\;\;\mbox{and}\;\; f(x,y)< c\}= \textup{proj}_X\bigl(D\cap f^{-1}([-\infty,c))\bigr).
\end{eqnarray*}
Proposition \ref{base_hierarchy_pj} (iii) and \eqref{Delta_pj} imply that 
$$D\cap f^{-1}([-\infty,c)) \in \Delta_q^1(X\times Y) \subset  \Sigma_q^1(X\times Y).$$
Thus, using (i) in Proposition \ref{base_hierarchy_pj}  applied to $\textup{proj}_X$, which is Borel measurable, and then \eqref{eq(v)_pj}, we get that 
 $$f_*^{-1}([-\infty,c))= \textup{proj}_X\bigl(D\cap f^{-1}([-\infty,c))\bigr)\in \Sigma_{q}^1(X)\subset \Delta_{q+1}^1(X).$$
  As $p$ and $q$ do not depend on $c$, Remark \ref{remarkxx} shows that $f_*$ is projective. The proof for $f^*$ is similar since
$$(f^*)^{-1}((c,+\infty])=\textup{proj}_X\bigl(D\cap f^{-1}((c,+\infty])\bigr).\Halmos$$
 \endproof
\section{Axiom of Projective Determinacy}
\label{pd_exp_lab}
We present in more detail the axiom of Projective Determinacy (PD), which will allow us to prove further results on projective functions.  Let $A \subset \mathcal{N}=\mathbb{N}^{\mathbb{N}}$ be a non-empty set. Imagine a two-player infinite game played as follows :
\begin{eqnarray*}
\mbox{I}\quad a_0 \quad \quad a_1 \quad \quad \cdots \label{game_pj}\\ 
\mbox{II} \quad \quad b_0 \quad \quad b_1 \quad \cdots\nonumber
\end{eqnarray*}
Player I plays $a_0\in \mathbb{N}$, then Player II plays $b_0\in \mathbb{N}$, then Player I plays $a_1\in \mathbb{N}$, etc. A play is thus a sequence $(a_0,b_0,a_1,b_1,\cdots) \in \mathcal{N}$. 
We say that Player I wins the game if  $(a_0,b_0,a_1,b_1,\cdots) \in A$. Else, if $(a_0,b_0,a_1,b_1,\cdots) \in  \mathcal{N} \setminus A,$ 
Player II wins. A strategy for Player I is a function $\sigma$ taking values in $\mathbb{N}$ and defined on the set of even-length integer sequences. 
For the sequence $(a_0,b_0,a_1,b_1,\cdots)$, the strategy $\sigma$ of Player I is defined as follows: $a_0 = \sigma(\emptyset)$, $a_1 = \sigma((a_0,b_0))$, $a_2 = \sigma((a_0,b_0,a_1,b_1))$, and so on. Note that if Player I follows the strategy $\sigma$, then the complete play is described by both $\sigma$ and $b=(b_0,b_1,\cdots),$ the actions of Player II, and is denoted by $\sigma * b$. A winning strategy $\sigma$  for Player I is a strategy under which Player I always wins, i.e., $\sigma * b \in A$ for all $ b\in\mathcal{N}$. Similarly, a strategy for Player II is a function $\tau$ taking values in $\mathbb{N}$ and defined 
on the set of odd-length integer sequences.  By $a * \tau$, we denote the play where Player I plays $a$ and Player II plays according to strategy $\tau$. A winning strategy for Player II is a strategy $\tau$ under which Player II always wins, i.e., $a * \tau
\in  \mathcal{N} \setminus A$ for all $a\in\mathcal{N}$. \\
 We say that $A$ is \textit{determined} if a winning strategy exists for one of the two players. 
One may wonder what kind of sets are determined and if a set exists that is not determined.  It turns out that any answer to these questions depends on the axiomatic set theory we consider. In the usual ZFC, any closed set $A$ is determined (see \citep[Corollary 7]{refgale_pj}) 
and, more generally, any Borel set is determined (see \citep[Theorem, p371]{refboreldet_pj}). However, as recalled in the introduction, the works of G{\"o}del and Solovay show that the determinacy of Borel sets is the best possible result provable in ZFC. 
For the second question, there exists a set that is not determined in ZFC,  see \citep[Lemma 33.1, p628]{refjech_pj}. 
\begin{remark}
Imagine a finite version of the two-player game defined above with $2N+2$ steps (instead of infinity). Let $A \subset \mathbb{N}^{2N+2}$. Then, ``Player I has a winning strategy" if and only if 
\begin{eqnarray}
\exists a_0\in \mathbb{N}, \forall b_0  \in \mathbb{N}, \exists a_1\in \mathbb{N}, \forall   b_1  \in \mathbb{N}, \cdots , \exists a_N 
 \in \mathbb{N}, \forall b_N  \in \mathbb{N}, \; (a_0,b_0,a_1,b_1,\cdots, a_N,b_N) \in A. \label{temp_finite_det}
\end{eqnarray}
Taking the contraposition of Assertion \eqref{temp_finite_det}, Morgan\rq{}s law shows that ``Player I has no winning strategy" if and only if 
\begin{eqnarray}
\forall a_0 \in \mathbb{N}, \exists b_0  \in \mathbb{N}, \forall a_1  \in \mathbb{N}, \exists   b_1  \in \mathbb{N}, \cdots , \forall a_N  \in \mathbb{N}, \exists b_N  \in \mathbb{N}, \; (a_0,b_0,a_1,b_1,\cdots, a_N,b_N) \in \mathbb{N}^{2N+2} \setminus A. \label{temp_finite_det2}
\end{eqnarray}
But, Assertion \eqref{temp_finite_det2} is exactly ``Player II has a winning strategy". 
Thus, the set $A$ is always determined, 
and all finite sets $A$ are determined in a finite sense. 
In the infinite case, $A$ is determined among to say that Morgan's law applies (formally) to countable sequences of $\forall$ and $\exists$. \\
\label{finit_det_pj}
\end{remark}

\begin{axiom}
Let $A \in \textbf{P}(\mathcal{N}).$ Then, $A$ is determined.\\
\label{projective_determinacy_pj}
\end{axiom}
One may wonder if adding the (PD) axiom to the other axioms of ZFC is legitimate. As a first observation, the (PD) axiom cannot be proved in ZFC. Indeed, otherwise a projective set will be universally measurable (see Theorem \ref{proj_is_univ_pj} below) and this is not true in ZFC. Indeed,  there exists a $\Sigma_2^1$ set, in the model $L$ of constructible sets,  that is not Lebesgue measurable. This was announced by G{\"o}del in 1938 and proved later by P. Novikov and J.W. Addison (see \citep[Corollary 25.28,  p495]{refjech_pj}). \\
The (PD) axiom is fruitful and extends properties of analytic sets holding in ZFC to projective sets. 
Moreover, every uncountable projective set contains a Cantor set (and so satisfies the Continuum hypothesis) and 
every projective set has the Baire property (see \citep[Theorem 38.17, p326]{refproj_pj}). 
What makes the (PD) axiom \lq\lq{}plausible\rq\rq{} is that it is implied by numerous set-theoretical statements, as the Proper Forcing Axiom,  and many of them come from areas of set theory with apparently no connection with projective sets. Examples of such statements can be found in \citep{refwood1_pj}. A remarkable discovery was the close link between the (PD) axiom and the notion of large cardinals: 
\citep{reflargedet_pj} proves that the existence of infinitely many Woodin cardinals implies the (PD) axiom. Conversely, assuming the (PD) axiom, one can construct inner models, which exhibit sequences of Woodin cardinals (see \citep[p573]{refwood1_pj}).

\section{Consequences  of the (PD) axiom}
\label{subsec_csq_pd_pj}
We now present two crucial consequences of the (PD) axiom (see Axiom \ref{projective_determinacy_pj}). The first one asserts that projective sets are universally measurable, and the second one that the class of projective sets satisfies the uniformisation property. 
This is indeed to get these two properties that the (PD) axiom was introduced. 
One may wonder if a kind of reciprocal holds, i.e., if we assume that all projective sets are Lebesgue measurable, have the Baire property, and can be uniformised by projective functions, does the (PD) axiom hold? This question, asked by Woodin in 1981,  and known as the 12th Delfino problem, was solved in 1997 by Steel, and the answer is no. 
So, in the rest of the paper, we could, instead of the (PD) axiom, assume that the projective sets are universally measurable and that $\Pi_{2n+1}^1(X \times Y)$ has the projective uniformisation property. This is implied, for example, by the existence of an infinite sequence $\kappa_0 <\kappa_1<\ldots$ of cardinals, with supremum $\lambda$, such that each $\kappa_n$ is $\lambda$-strong. This will be a weaker assumption, but as explained above, the (PD) axiom is a plausible and well-accepted axiom.

We define the universal sigma-algebra on $X$ as
\[
\mathcal{B}_c(X) := \bigcap_{P \in \mathfrak{P}(X)} \mathcal{B}_P(X),
\]
where $\mathfrak{P}(X)$ is the set of probability measures on $X$ and $\mathcal{B}_P(X)$ denotes the completion of $\mathcal{B}(X)$ with respect to $P \in \mathfrak{P}(X)$. 
\begin{theorem}
Assume the (PD) axiom. Let $X$ be a Polish space. Then, $\textbf{P}(X) \subset \mathcal{B}_c(X)$.
\label{proj_is_univ_pj}
\end{theorem}
\proof{Proof.}
See \citep[Remark 2, p71]{refprojareuniv1_pj} and \citep[Theorem 38.17, p326]{refproj_pj}.\Halmos \\ \endproof
An accessible proof of Theorem \ref{proj_is_univ_pj} for $\Sigma_2^1(X)$ sets is available in \citep[p308]{refproj_pj}. Interestingly, \citep{refmai_pj} shows that the universal measurability of all value functions associated with a Borel gambling problem is equivalent to the sets of $\Sigma_2^1(X)$ being universally measurable.\\
Theorem \ref{proj_is_univ_pj} implies that  $\Sigma_n^1(X)$, $\Pi_n^1(X)$ and $\Delta_n^1(X)$ are included in $\mathcal{B}_c(X)$ for all $n\geq 1$, see \eqref{eq(v)_pj} in Proposition \ref{base_hierarchy_pj}. It also implies that projective functions are universally measurable. 
\begin{proposition}
Assume the (PD) axiom. 
Let $X$ and $Y$ be Polish spaces. If $f:X \to Y$ is a projective function, then $f$ is universally measurable (see \citep[Definition 7.20, p171]{ref1_pj}).
\label{proj_mes_pj}
\end{proposition}
\proof{Proof.}
This follows directly from Theorem \ref{proj_is_univ_pj} as for all $n\geq 1$, $\Delta_n^1(X) \subset \textbf{P}(X)\subset \mathcal{B}_c(X)$.\Halmos \\ \endproof
We now give the second consequence of the (PD) axiom that we will use. This is the so-called uniformisation theorem due to Moschovakis,  see  \cite{refunifofproj_pj}. 
\begin{theorem}
Assume the (PD) axiom. Let $X$ and $Y$ be Polish spaces, $n\geq 0$ and $A \in \Pi_{2n+1}^1(X \times Y)$. Then, there exists $A^* \in \Pi_{2n+1}^1(X \times Y)$ such that $A^* \subset A$ and for all $x\in X$,
\begin{eqnarray}
\exists y \in Y,\; (x,y)\in A \iff \exists ! y\in Y,\; (x,y)\in A^*.
\label{unif_eq_pj}
\end{eqnarray}
The set $A^*$ is called a $\Pi_{2n+1}^1(X \times Y)$-uniformisation of $A$. 
\label{unif_thm_brut_pj}
\end{theorem}
\proof{Proof.}
See \citep[p120]{refproj_pj} for the definition of a uniformisation and \citep[Theorem 1]{refunifofproj_pj},  \citep[Theorem 6C.5]{Moschovakis} or \citep[Corollary 39.9, p339]{refproj_pj} for the proof.
\Halmos \\\endproof
\begin{remark} 
Theorem \ref{unif_thm_brut_pj} provides a uniformisation of any set $A \in \Pi_{2n+1}^1(X\times Y)$. 
The proof of Theorem \ref{unif_thm_brut_pj} is much more complex than that of Theorem \ref{proj_is_univ_pj}. It relies on the notion of scale. A scale on a given set $A$ is a sequence of ordinal functions on $A$ that satisfy sufficient continuity conditions to select an element out of $A$. The reader can consult \citep[36(B),  p299]{refproj_pj} for a complete description of scales. Note that when $n=0$, i.e., for the class $\Pi_1^1(X\times Y)$ of coanalytic sets and also the class $\Sigma_2^1(X\times Y)$, Theorem \ref{unif_thm_brut_pj} can be shown inside ZFC (see the \cite[Novikov-Kondo-Addison Uniformization Theorem 4E.4]{Moschovakis}). \\
\label{rem_uni_th_pj}
\end{remark}
Theorem \ref{unif_thm_brut_pj} implies directly that, under the (PD) axiom, a projective selection on a projective set is always possible. This proposition can be seen as an extension of the Jankov von Neuman theorem \citep[Proposition 7.49, p182]{ref1_pj}, which holds for analytic sets in ZFC.
\begin{proposition}
Assume the (PD) axiom. 
Let $X$ and $Y$ be Polish spaces. 
 Let $A \in \mathbf{P}(X \times Y)$. Then, there exists a projective function $\phi : \textup{proj}_X(A) \to Y$ such that $\textup{Graph}(\phi) \subset A$.
\label{univ_sel_proj_pj}
\end{proposition}
\proof{Proof.}
Let $A \in \mathbf{P}(X \times Y)$. Using \eqref{Delta_pj} and (iv) in Proposition \ref{base_hierarchy_pj},  there exists some $n\geq 0$ such that 
$$A\in \Delta_n^1(X\times Y)\subset \Pi_n^1(X\times Y) \subset \Pi_{2n+1}^1(X\times Y).$$ 
Applying now Theorem \ref{unif_thm_brut_pj} to $A$, there exists $A^* \in \Pi_{2n+1}^1(X\times Y)$ such that $A^*\subset A$ and \eqref{unif_eq_pj} holds true. So, we get that
\begin{eqnarray*}
x\in\textup{proj}_X(A) &\iff& \exists y_0\in Y\; \mbox{such that } (x,y_0)\in A  \iff 
 \exists !  y\in Y\; \mbox{such that } (x,y)\in A^*.
\end{eqnarray*}
For all $x\in\textup{proj}_X(A)$, we set $\phi(x) := y$. Then, $\phi$ is a function because of the unicity obtained in \eqref{unif_eq_pj}. By definition of $\phi$, 
$$\textup{Graph}(\phi) := \{(x,y)\in X\times Y,\; x\in\textup{proj}_X(A),\; y=\phi(x)\}\subset A^*\subset A.$$ 
As $A \in \mathbf{P}(X \times Y)$, (vi) in Proposition \ref{base_hierarchy_pj} applied to $\textup{proj}_X$ shows that $\textup{proj}_X(A)\in \mathbf{P}(X).$
As $\textup{Graph}(\phi)\in \mathbf{P}(X \times Y)$, (iii) in Proposition \ref{pareil} shows that $\phi$ is a projective function. \Halmos \\ \endproof

\section{Properties of projective functions under the (PD) axiom}
\label{subsec_csq_pd2_pj}
This section provides two results for projective functions that necessitate the (PD) axiom 
(or to assume that the projective sets are universally measurable and that $\Pi_{2n+1}^1(X \times Y)$ has the projective uniformisation property). We show that projectively measurable $\epsilon$-optimal selectors exist for projective functions  (see Proposition \ref{prop_sel_mes_proj_pj}). Then, we prove that the integral of a projective function remains a projective function under the $(-\infty)$ integration (see Proposition \ref{prj_cvt_pj}), and that the multiple-priors conditional support of a projective function is a projective set-valued mapping (see Proposition \ref{aff_is_proj_pj}).

\subsection{Existence of $\epsilon$-optimizer for projective functions}
Assuming the (PD) axiom, we show that $\epsilon$-optimal selectors exist for the optimisation of projective functions. This proposition can be seen as an extension of  \citep[Proposition 7.50, p184]{ref1_pj}, where \eqref{mes_sel_inf_pj} (resp.  \eqref{mes_sel_sup_pj}) holds for lower (resp. upper)-semianalytic functions.

\begin{proposition}
Assume the (PD) axiom. Let $X$ and $Y$ be Polish spaces. Let $D \in \textbf{P}(X\times Y)$ and $f : X \times Y \to \mathbb{R}\cup\{-\infty,+\infty\}$ be a projective function.  Let $\epsilon_*>0$ and $\epsilon^*>0$. Then, there exist projective functions $\phi_*, \phi^* : \textup{proj}_X(D) \to Y$ such that $\textup{Graph}(\phi_*) \subset D$, $\textup{Graph}(\phi^*) \subset D$ and for all $x\in \textup{proj}_X(D)$,
\begin{eqnarray}
\label{mes_sel_inf_pj}
f(x,\phi_*(x)) & < &\left\{
   \begin{array}{ll}
  f_*(x) +\epsilon_* \;\;\mbox{if}\;\; f_*(x)>-\infty, \\
  -\frac{1}{\epsilon_*} \;\;\mbox{if}\;\; f_*(x)=-\infty\\
    \end{array}
\right.\\
\label{mes_sel_sup_pj}
f(x,\phi^*(x)) & > &\left\{
   \begin{array}{ll}
  f^*(x) -\epsilon^* \;\;\mbox{if}\;\; f^*(x)<+\infty, \\
  \frac{1}{\epsilon^*} \;\;\mbox{if}\;\; f^*(x)=+\infty,\\
    \end{array}
\right.
\end{eqnarray}
where  $D_x := \{y\in Y,\; (x,y)\in D\}$ for $x\in X$ and $f_*, f^* : \textup{proj}_X(D) \to \mathbb{R}\cup\{-\infty,+\infty\}$ are defined by
\begin{eqnarray*}
f_*(x) := \inf_{y\in D_x} f(x,y) \;\; \mbox{and} \;\; f^*(x) := \sup_{y\in D_x} f(x,y).
\end{eqnarray*}
\label{prop_sel_mes_proj_pj}
\end{proposition}

\proof{Proof.}
We show \eqref{mes_sel_inf_pj}. Let $\epsilon_*>0$ and $E:= (A_1 \cap B_1) \cup (A_2\cap B_2) \subset X\times Y$ where,
\begin{eqnarray*}
A_1&:=& \left\{(x,y)\in D,\; f(x,y)<f_*(x)+\epsilon_*\right\}\\
A_2&:=& \left\{(x,y)\in D,\; f(x,y)<-1/\epsilon_*\right\}\\
B_1 &:=& \{x\in \textup{proj}_X(D),\; f_*(x)>-\infty\}\times Y\\
B_2 &:=& \{x\in \textup{proj}_X(D),\; f_*(x)=-\infty\}\times Y.
\end{eqnarray*}
We prove first that $$\textup{proj}_X(E) = \textup{proj}_X(D).$$ The ``$\subset$" inclusion is trivial as $E\subset D$. Now, let $x\in \textup{proj}_X(D).$ 
Then, $D_x \neq \emptyset$. If $f_*(x)>-\infty$, by definition of the infimum, there exists  $y\in D_x$ such that $f(x,y)<f_*(x)+\epsilon_*$ and  $(x,y)\in A_1 \cap B_1 \subset E$. Similarly, if $f_*(x)=-\infty$, there exists $y\in D_x$ such that $(x,y)\in A_2 \cap B_2 \subset E$. Thus, $x\in \textup{proj}_X(E)$ and the reverse inclusion is proved. \\
Assume for a moment that $E \in \textbf{P}(X\times Y)$. Let $\phi_* : \textup{proj}_X(E) \to Y$ be the function such that $\textup{Graph}(\phi_*)\subset E$ given by Proposition \ref{univ_sel_proj_pj}. Then, \eqref{mes_sel_inf_pj} holds true for $\phi_*$ for all $x\in \textup{proj}_X(E)=\textup{proj}_X(D)$ and $\phi_*$ is the desired function as $E\subset D$. \\
It remains to show that $E \in \textbf{P}(X\times Y)$. As $f$ is projective and $D \in \textbf{P}(X\times Y)$, Proposition \ref{base_hierarchy_pj} (vi) shows that 
$$A_2 = f^{-1}([-\infty,-1/\epsilon_*))\cap D\in \textbf{P}(X\times Y).$$
As $f_*$ is a projective function (see Proposition \ref{prop_infsup_mes_proj_pj}), Proposition \ref{base_hierarchy_pj} (vi) and Proposition \ref{lemma_real_proj_pj} show that 
$(x,y) \mapsto f(x,y)-f_*(x)$ is also projectively measurable: $A_1 \in \textbf{P}(X\times Y)$. 
 Now, Proposition \ref{base_hierarchy_pj} (vi) again shows that $\textup{proj}_X(D)$,  $\{x\in \textup{proj}_X(D),\; f_*(x)>-\infty\}$ and $\{x\in \textup{proj}_X(D),\; f_*(x)=-\infty\}$ belong to $\textbf{P}(X)$,  that $B_1$ and $B_2$ belong to $\textbf{P}(X \times Y)$ and also that $E \in \textbf{P}(X\times Y)$. \\
 Finally, as $f^* = -(-f)_*$, \eqref{mes_sel_sup_pj} follows from \eqref{mes_sel_inf_pj} applied to $-f$, which is projectively measurable, see Proposition \ref{lemma_real_proj_pj}.
\Halmos \endproof

\subsection{Integration of projective functions}
\label{inf_cvt_sub_pj}
Under the (PD) axiom, any projective set $A$ is universally measurable (see Theorem \ref{proj_is_univ_pj}). So, $p(A)$ can be defined for any probability measure $p$ and, more generally, the classical measure theory results can be applied in the projective setup. Indeed, any projective function $f$ is universally measurable (see Proposition \ref{proj_mes_pj}). So, $\int f dp$ is well-defined and Fubini's theorem (see \citep[Proposition 7.45 p175]{ref1_pj})  can be used for projective functions and projectively measurable stochastic kernels. 
The main result of this section is an extension of \citep[Proposition 7.43, p169]{ref1_pj}  to the projective hierarchy (see Proposition \ref{prop_crutial_hierarchy_pj}). 
This result is crucial to show that integrals of projective functions remain projective (see Proposition \ref{prj_cvt_pj}), and also that the multiple-priors conditional support of a projective function is a projective set-valued mapping (see Proposition \ref{aff_is_proj_pj}). To show that, we first establish a lemma which roughly states that the sets in $\Sigma_{n+1}^1(X)$, and also the projective sets (see \eqref{eq(v)_pj}), are close to analytic sets in terms of measure. 
\begin{lemma}
Assume the (PD) axiom. Let $X$ be a Polish space, $\mu\in \mathfrak{P}(X)$,  $n\geq 1$ and $A := \textup{proj}_X(C)$ for some $C\in\Pi_{n}^1(X\times \mathcal{N})$. Then, there exists $C' \in \mathcal{B}(X\times \mathcal{N})$ such that $C'\subset C$ and $A \setminus \textup{proj}_X(C')$ is a $\mu$-null-set.
\label{lemma_before_crit_pj}
\end{lemma}
\proof{Proof.} Recalling \eqref{Sigma_pj}, \eqref{eq(v)_pj} and Theorem \ref{proj_is_univ_pj}, 
$$A \in \Sigma_{n+1}^1(X)\subset \Delta_{n+2}^1(X)\subset \textbf{P}(X) \subset \mathcal{B}_c(X)$$
and  there exists some $A' \in \mathcal{B}(X)$ such that $A' \subset A$ and $A \setminus A'$ is a $\mu$-null-set. \\
As $C\in \Pi_{n}^1(X\times \mathcal{N}) \subset \Delta_{n+1}^1(X\times \mathcal{N}) \subset \textbf{P}(X\times \mathcal{N})$ (see \eqref{eq(v)_pj}), Proposition \ref{univ_sel_proj_pj} provides a projective function $\phi : A=\textup{proj}_X(C)  \to \mathcal{N}$ such that $\textup{Graph}(\phi)\subset C$. 
Proposition \ref{proj_mes_pj} shows that $\phi$ is also universally measurable.  
Let 
$$\widetilde{\phi}(x):=\phi(x) \mbox{ if }x\in A'\mbox{ and }\widetilde{\phi}(x) := u_0 \mbox{ if }x\in X\setminus A',$$ where $u_0$ is an arbitrary element of $\mathcal{N}$. Then, $\widetilde{\phi}$ is also universally measurable. Thus, \citep[Lemma 7.27, p173]{ref1_pj}  provides a $\mathcal{B}(X)$-measurable function $\widetilde{\psi} : X \to \mathcal{N}$ such that $\widetilde{\phi} =\widetilde{\psi}$ $\mu$-a.s.  Let $B' \in \mathcal{B}(X)$ such that $\mu(X \setminus B')=0$  and $\widetilde{\phi}(x) =\widetilde{\psi}(x)$ for all $x\in B'$. Let 
$$C' := \{(x,u)\in X\times \mathcal{N},\; u= \widetilde{\psi}(x)\;\; \mbox{and} \;\; x\in A' \cap B'\}.$$ 
Then, $C' \in \mathcal{B}(X\times \mathcal{N})$ (see \citep[Corollary 7.14.1, p121]{ref1_pj}). Moreover, $C'\subset C$. Indeed, let $(x,u)\in C'$. Then, $u= \widetilde{\psi}(x) = \widetilde{\phi}(x) =\phi(x)$ as $x\in A' \cap B'$. So, $(x,u)\in \textup{Graph}(\phi)\subset C$ and $C'\subset C$. \\
Let $x\in A\setminus\textup{proj}_X(C')$. If $x\in A' \cap B'$, then $\phi(x)=\widetilde{\phi}(x) = \widetilde{\psi}(x)$ and $x\in \textup{proj}_X(C')$ which is a contradiction. Thus, 
$$A\setminus \textup{proj}_X(C') \subset A \cap (X\setminus (A' \cap B'))\subset (A \setminus A') \cup (X \setminus B')$$ and $A\setminus\textup{proj}_X(C')$ is indeed a $\mu$-null-set.
\Halmos\\ \endproof

The next proposition extends \citep[Proposition 7.43, p169]{ref1_pj}, which holds for $\Sigma_1^1(X)$, i.e., analytic sets,  to $\Sigma_n^1(X)$ and thus projective sets. The proof completely differs from that of \citep[Proposition 7.43, p169]{ref1_pj} and is based on Lemma \ref{lemma_before_crit_pj} as well as on a parametrisation of $\Delta_1^1(X \times \mathcal{N})$ given in \cite{refproj_pj}.
\begin{proposition}
Assume the (PD) axiom. Let $X$ be a Polish space, $n\geq 1$ and $A \in \Sigma_n^1(X)$. Then, for all $r\in \mathbb{R}$,
\begin{eqnarray}
W_r := \{\mu\in \mathfrak{P}(X),\; \mu[A]\geq r\}\in \Sigma_n^1(\mathfrak{P}(X)).
\label{crucial_mesure_hierarchy_pj}
\end{eqnarray}
\label{prop_crutial_hierarchy_pj}
\end{proposition}
\proof{Proof.}
The case $n=1$ is shown  in \citep[Proposition 7.43, p169]{ref1_pj}.  We now prove \eqref{crucial_mesure_hierarchy_pj} for  $A\in\Sigma_{n+1}^1(X)$ with  $n\geq 1$. There exists $C\in\Pi_n^1(X\times \mathcal{N})$ such that $A=\textup{proj}_X(C)$. We have that
\begin{eqnarray}
W_r = \{\mu\in\mathfrak{P}(X),\; \exists C' \in \mathcal{B}(X\times \mathcal{N}),\; C'\subset C \;\; \mbox{and} \;\;  \mu[\textup{proj}_X(C')]\geq r\}.
\end{eqnarray}
Indeed, the ``$\subset$" inclusion is a direct consequence of Lemma \ref{lemma_before_crit_pj}, while  the ``$\supset$" one is immediate from the monotony of $\mu$.\\ 
We now introduce a parametrisation of $\mathcal{B}(X\times \mathcal{N})=\Delta_1^1(X\times \mathcal{N})$ (see (iii) in Proposition \ref{base_hierarchy_pj}). Let $\mathcal{C}:= \{0,1\}^{\mathbb{N}}$ be the Cantor space (see \citep[p13]{refproj_pj}). Using \citep[Theorem 35.5, p283]{refproj_pj}, there exist $\mathcal{D}\in \Pi_1^1(\mathcal{C})$ and $S\in \Sigma_1^1(\mathcal{C} \times X\times \mathcal{N})$ such that $\Delta_1^1(X\times \mathcal{N}) = \{S_d,\; d\in \mathcal{D}\}$ where $S_d := \{(x,u)\in X\times \mathcal{N},\; (d,x,u)\in S\}$ for all $d\in \mathcal{C}$. Thus,
$$W_r = \{\mu\in\mathfrak{P}(X),\; \exists d\in \mathcal{D},\; S_d \subset C \;\; \mbox{and} \;\;  \mu[\textup{proj}_X(S_d)] \geq r\}=\textup{proj}_{\mathfrak{P}(X)}(\mathcal{Z}_r^1 \cap \mathcal{Z}_r^2 \cap \mathcal{Z}_r^3),$$
where
\begin{eqnarray*}
\mathcal{Z}_r^1 &:=& \{(\mu,d)\in\mathfrak{P}(X)\times \mathcal{C},\;d\in \mathcal{D}\}\\
\mathcal{Z}_r^2 &:=&\{(\mu,d)\in\mathfrak{P}(X)\times \mathcal{C},\;S_d \subset C\}\\
\mathcal{Z}_r^3 &:=&\{(\mu,d)\in\mathfrak{P}(X)\times \mathcal{C},\; \mu[\textup{proj}_X(S_d)]\geq r\}.
\end{eqnarray*}
Assume for a moment that $\mathcal{Z}_r^1$, $\mathcal{Z}_r^2$ and $\mathcal{Z}_r^3$ belong to $\Sigma_{n+1}^1(\mathfrak{P}(X)\times \mathcal{C})$. Then (i) in Proposition \ref{base_hierarchy_pj} applied to $\textup{proj}_{\mathfrak{P}(X)}$, which is Borel measurable, ensures that $W_r\in \Sigma_{n+1}^1(\mathfrak{P}(X)),$ which is the desired result.

\noindent \textit{(i) $\mathcal{Z}^1_r \in \Sigma_{n+1}^1(\mathfrak{P}(X)\times \mathcal{C}).$} \\
(v) in Proposition \ref{base_hierarchy_pj} implies that 
 $$\mathcal{Z}^1_r = \mathfrak{P}(X) \times \mathcal{D} \in \Pi_1^1(\mathfrak{P}(X))\times \Pi_1^1(\mathcal{C}) \subset \Pi_1^1(\mathfrak{P}(X)\times \mathcal{C})$$ and we conclude as $\Pi_1^1 \subset \Pi_n^1 \subset \Delta_{n+1}^1 \subset \Sigma_{n+1}^1$, see (iv) and \eqref{eq(v)_pj} in Proposition \ref{base_hierarchy_pj}.

\noindent \textit{(ii) $\mathcal{Z}_r^2 \in \Sigma_{n+1}^1(\mathfrak{P}(X)\times \mathcal{C}).$} \\
Taking the complement of $\mathcal{Z}_r^2$, we obtain that
\begin{eqnarray*}
(\mathfrak{P}(X)\times \mathcal{C})\setminus \mathcal{Z}_r^2 &=& \{(\mu,d)\in\mathfrak{P}(X)\times \mathcal{C},\; \exists (x,u)\in X\times \mathcal{N},\; (d,x,u) \in S \;\; \mbox{and} \;\; (x,u)\in (X\times \mathcal{N}) \setminus C\}\\
&=& \textup{proj}_{\mathfrak{P}(X)\times \mathcal{C}}\big(\mathfrak{P}(X) \times S\big) \cap \Bigl(\mathfrak{P}(X) \times \mathcal{C} \times \bigl((X\times \mathcal{N}) \setminus C)\bigr)\Bigr).
\end{eqnarray*}
As $S\in \Sigma_1^1(\mathcal{C}\times X\times \mathcal{N}) \subset \Sigma_n^1(\mathcal{C}\times X\times \mathcal{N})$ (see (iv) in Proposition \ref{base_hierarchy_pj}),   (v) in Proposition \ref{base_hierarchy_pj} implies that  
$$\mathfrak{P}(X) \times S \in \Sigma_n^1(\mathfrak{P}(X) \times \mathcal{C} \times X \times \mathcal{N}).$$ Now, as $C\in \Pi_n^1(X\times \mathcal{N})$, \eqref{Pi_pj} shows that $(X\times \mathcal{N})\setminus C \in \Sigma_n^1(X\times \mathcal{N})$. So, we deduce from (v) and (i) in Proposition \ref{base_hierarchy_pj} that 
\begin{eqnarray*}
\mathfrak{P}(X) \times \mathcal{C} \times \bigl((X \times \mathcal{N}) \setminus C\bigr) & \in &  \Sigma_n^1(\mathfrak{P}(X) \times \mathcal{C} \times X \times \mathcal{N}) \\
\bigl(\mathfrak{P}(X) \times S\bigr) \cap \Bigl(\mathfrak{P}(X) \times \mathcal{C} \times \bigl((X\times \mathcal{N}) \setminus C)\bigr)\Bigr) & \in &  \Sigma_n^1(\mathfrak{P}(X) \times \mathcal{C} \times X \times \mathcal{N})\\
(\mathfrak{P}(X)\times \mathcal{C}) \setminus \mathcal{Z}_r^2 =\textup{proj}_{\mathfrak{P}(X)\times \mathcal{C}}\bigl(\mathfrak{P}(X) \times S\bigr) \cap \Bigl(\bigl(\mathfrak{P}(X) \times \mathcal{C} \times ((X\times \mathcal{N}) \setminus C)\bigr)\Bigr) &\in &  \Sigma_n^1(\mathfrak{P}(X)\times \mathcal{C}).
\end{eqnarray*} 
Thus, using \eqref{Pi_pj} and \eqref{eq(v)_pj}, $\mathcal{Z}_r^2 \in \Pi_n^1(\mathfrak{P}(X)\times \mathcal{C}) \subset \Sigma_{n+1}^1(\mathfrak{P}(X)\times \mathcal{C})$.

\noindent \textit{(iii) $Z_r^3 \in \Sigma_{n+1}^1(\mathfrak{P}(X)\times \mathcal{C}).$}\\
Remark first that for all $d\in\mathcal{C}$, 
$$\textup{proj}_X(S_d) = \bigl(\textup{proj}_{\mathcal{C}\times X}(S)\bigr)_d.$$ 
As $S\in\Sigma_1^1(\mathcal{C}\times X \times \mathcal{N})$, (i) in Proposition \ref{base_hierarchy_pj} shows that $\textup{proj}_{\mathcal{C}\times X}(S) \in \Sigma_1^1(\mathcal{C}\times X)$. Note now that for all $(\mu,d) \in \mathfrak{P}(X)\times\mathcal{C}$,  using Fubini's theorem, we obtain that 
$$\mu[\textup{proj}_X(S_d)] = \mu\bigl[(\textup{proj}_{\mathcal{C}\times X}(S))_d \bigr]= (\delta_d \otimes \mu)[\textup{proj}_{\mathcal{C}\times X}(S)],$$
where $\delta_d$ is the Dirac measure in $d$, i.e., $\delta_d[A]=1$ if $d\in A$ and 0 otherwise, for any $A \subset \mathcal{C}.$ 
So, 
$$\mathcal{Z}_r^3 = \zeta^{-1}\bigl[\{\nu\in\mathfrak{P}(\mathcal{C}\times X),\; \nu[\textup{proj}_{\mathcal{C}\times X}(S)] \geq r\}\bigr],$$
where $\zeta : \mathfrak{P}(X) \times \mathcal{C} \to \mathfrak{P}(\mathcal{C} \times X)$ is defined by $\zeta(\mu,d):=\delta_d \otimes \mu$ for all $(\mu,d)\in \mathfrak{P}(X) \times \mathcal{C}$. Then, $\zeta$ is continuous (see \citep [Lemma 7.12, p144]{ref1_pj} and \citep [Corollary 7.21.1, p130]{ref1_pj}). Recalling that $\textup{proj}_{\mathcal{C}\times X}(S) \in \Sigma_1^1(\mathcal{C}\times X)$ and using \eqref{crucial_mesure_hierarchy_pj} for $n=1$, we obtain that 
$$\{\nu\in\mathfrak{P}(\mathcal{C}\times X),\; \nu[\textup{proj}_{\mathcal{C}\times X}(S)] \geq r\} \in \Sigma_1^1(\mathfrak{P}(\mathcal{C}\times X)).$$ 
Thus, $\mathcal{Z}_r^3 \in \Sigma_1^1(\mathfrak{P}(X) \times \mathcal{C}) \subset \Sigma_{n+1}^1(\mathfrak{P}(X) \times \mathcal{C})$ using (i) and (iv) in Proposition \ref{base_hierarchy_pj}.
\Halmos \\ \endproof
We now show that the integral of a projective function with respect to a projectively measurable stochastic kernel is a projective function. Proposition \ref{prj_cvt_pj} 
is an extension of \citep[Proposition 7.48, p180]{ref1_pj}, which holds for a lower semi-analytic function and a Borel measurable stochastic kernel, under the reverse convention to \eqref{cvt_inf_pj}, i.e.,  
$+\infty-\infty =+\infty$.  Proposition \ref{prop_crutial_hierarchy_pj} will be in force. Here, we define generalised integrals as in \citep{cf2024}. Indeed, the integral $\int_{-}$ is adapted to the maximisation problems of Section \ref{secappli}, where the convention $+\infty-\infty =-\infty$  is usual. 
\begin{definition}
Let $X$ be a Polish space. Let $f : X \to \mathbb{R}\cup\{-\infty,+\infty\}$ be a  universally measurable function and let $p\in \mathfrak{P}(X)$. 
The $(-\infty)$ integral denoted by $\int_{-} f dp$ and the $(+\infty)$ integral denoted by $\int^{-} f dp$ are defined as follows. When $\int f^+ dp<+\infty$ or $\int f^- dp<+\infty$, both integrals equal  the extended integral of $f$, i.e.
\begin{eqnarray}
\int_{-} f dp = \int^{-} f dp:=\int f^+ dp - \int f^- dp. \label{base_integral_pj}
\end{eqnarray}
Otherwise, $\int_{-} f dp := -\infty \;\;\mbox{and}\;\; \int^{-} f dp := +\infty.$
\end{definition}
We now define projectively measurable stochastic kernels (see also \citep[Definition 7.12, p134]{ref1_pj}).
\begin{definition}
Let $X$ and $Y$ be Polish spaces. A stochastic kernel on $Y$ given $X$ is a function $q(\cdot|\cdot) : \mathcal{B}(Y)\times X \to \mathbb{R}$ such that  $q[\cdot|x] \in \mathfrak{P}(Y)$ for all $x \in X$. \\
(i) Let $r\geq 1.$ If  $X \ni x \mapsto q[\cdot|x]  \in \mathfrak{P}(Y)$ is $\Delta_r^1(X)$-measurable, then $q$ is called a  $\Delta_r^1(X)$-measurable stochastic kernel on $Y$ given $X$.\\
(ii) If $X \ni x \mapsto q[\cdot|x] \in \mathfrak{P}(Y)$ is a projective function, then $q$ is called a projectively measurable stochastic kernel on $Y$ given $X$.
\end{definition}

\begin{proposition}
Assume the (PD) axiom. Let $X$ and $Y$ be Polish spaces. Let $f : X \times Y \to \mathbb{R}\cup\{-\infty,+\infty\}$ and let $q$ be a stochastic kernel on $Y$ given $X$. Let $\lambda : X  \to \mathbb{R}\cup\{-\infty,+\infty\}$ be defined by $$\lambda(x):=\int_{-} f(x,y)q[dy|x].$$
\noindent (i) Assume that $q$ is a $\Delta_r^1(X)$-measurable stochastic kernel for some $r\geq 1$ and that $f$ is $\Delta_p(X\times Y)$-measurable for some $p\geq 1$. Then, $\lambda$ is $\Delta_{p+r+2}^1(X)$-measurable.\\
\noindent (ii) Assume that $q$ is a  projectively measurable stochastic kernel and that $f$ is a projective function. Then, $\lambda$ is a projective function.
\label{prj_cvt_pj}
\label{univ_cvt_pj} 
\end{proposition}
\proof{Proof.}
The proof of (i) is close to that of \citep[Proposition 7.48, p180]{ref1_pj}, but recall that \citep{ref1_pj} uses the reverse convention to \eqref{cvt_inf_pj}. Assume first that $f\geq 0$. In this case $\int_{-}f(x,y) q[dy|x] = \int f(x,y) q[dy|x]$. 
Let $\theta_f : \mathfrak{P}(X\times Y) \to \mathbb{R}\cup\{-\infty,+\infty\}$  and $\sigma : X \to \mathfrak{P}(X\times Y)$ be defined by 
$$\theta_f(\nu):= \int_{X\times Y} f d\nu \quad  \mbox{ and }  \quad \sigma(x):= \delta_{x} \otimes q[\cdot|x].$$  
Then,  $\lambda(x)= \theta_f(\sigma(x)).$ \\
First, we show that $\sigma$ is $\Delta_{r+1}^1(X)$-measurable. 
As $x \mapsto q[\cdot|x]$ is $\Delta_r^1(X)$-measurable and $x \mapsto \delta_x$ is continuous (see \citep[Corollary 7.21.1, p130]{ref1_pj}) and thus $\Delta_r^1(X)$-measurable (see (iii) and (iv) in Proposition \ref{base_hierarchy_pj}),  Lemma \ref{couple_fct_mes} (i) implies that 
$x \mapsto (\delta_{x},q[\cdot|x])$ is $\Delta_r^1(X)$-measurable. Recalling that $(\nu,\mu)\mapsto \nu \otimes \mu$ is continuous and thus $\mathcal{B}(\mathfrak{P}(X)\times \mathfrak{P}(Y)) = \Delta_1^1(\mathfrak{P}(X)\times \mathfrak{P}(Y))$-measurable (see (iii) in Proposition \ref{base_hierarchy_pj}), Proposition \ref{comp_proj_pj} (i)  proves that $\sigma$ is $\Delta_{r+1}^1(X)$-measurable. \\
Assume for a moment that $\theta_f$ is $\Delta_{p+1}^1(\mathfrak{P}(X\times Y))$-measurable. Then, $\lambda=\theta_f \circ \sigma$ is $\Delta_{p+r+2}^1(X)$-measurable using again Proposition \ref{comp_proj_pj} (i). 
We now prove that $\theta_f$ is $\Delta_{p+1}^1(\mathfrak{P}(X\times Y))$-measurable. 
For all $n\geq 0$, let 
$$f_n(x,y):=\min(n,f(x,y)) \mbox{ and } E_n := \{(x,y,b)\in X \times Y\times \mathbb{R},\; f_n(x,y)\leq b\leq n\}.$$ Then, $f_n \uparrow f.$  
Let $\mu$ be the Lebesgue measure on $\mathbb{R}$ and  let $\nu\in \mathfrak{P}(X\times Y).$ 
By Fubini's theorem,  we get that 
\begin{eqnarray*}
(\nu \otimes \mu)(E_n) &=& \int_{X\times Y} \int_{\mathbb{R}} 1_{E_n} d\mu d\nu =\int_{X\times Y} [n-f_n] d\nu = n- \int_{X\times Y} f_n d\nu.
\end{eqnarray*} 
Let  $c\in\mathbb{R}.$ As the sequence $(f_n)_{n\geq 0}$ is nondecreasing, the monotone convergence theorem shows that
\begin{eqnarray*}
\Bigl\{\nu\in \mathfrak{P}(X\times Y),\; \int_{X\times Y} f d\nu \leq c\Bigr\} &=& \bigcap_{n=0}^{+\infty} \Bigl\{\nu\in \mathfrak{P}(X\times Y),\; \int_{X\times Y} f_n d\nu \leq c\Bigr\}\\
&=& \bigcap_{n=0}^{+\infty}\{\nu\in \mathfrak{P}(X\times Y),\; (\nu \otimes \mu)[E_n] \geq n-c\}\\
&=& \bigcap_{n=0}^{+\infty} \Psi^{-1}\bigl(\{\rho\in \mathfrak{P}(X\times Y \times \mathbb{R}),\; \rho[E_n] \geq n-c\}\bigr),
\end{eqnarray*}
where $\Psi : \mathfrak{P}(X \times Y) \to \mathfrak{P}(X \times Y \times \mathbb{R})$ is such that $\Psi(\nu):=\nu \otimes \mu$ for all $\nu\in \mathfrak{P}(X \times Y)$. 
Fix $n\geq 0$. Assume for a moment that $E_n \in \Delta_p^1(X\times Y\times \mathbb{R}).$ Then, $E_n \in \Sigma_p^1(X\times Y\times \mathbb{R})$ and  Proposition \ref{prop_crutial_hierarchy_pj} shows that 
$$\{\rho\in\mathfrak{P}(X \times Y \times \mathbb{R}),\; \rho[E_n] \geq n-c \} \in \Sigma_p^1(\mathfrak{P}(X\times Y\times \mathbb{R})).$$ 
As $\Psi$ is continuous (see \citep[Lemma 7.12, p144]{ref1_pj}), (i) and \eqref{eq(v)_pj} in Proposition \ref{base_hierarchy_pj} show that 
$$\Bigl\{\nu\in \mathfrak{P}(X\times Y),\; \int_{X\times Y} f d\nu \leq c\Bigr\} \in \Sigma_p^1(\mathfrak{P}(X\times Y)) \subset \Delta_{p+1}^1(\mathfrak{P}(X\times Y))$$ 
and  $\theta_f$  is $\Delta_{p+1}^1(\mathfrak{P}(X\times Y))$-measurable (see Remark \ref{remarkxx}).\\
We still have to prove that $E_n \in \Delta_p^1(X\times Y\times \mathbb{R})$. Remark first that 
\begin{eqnarray*} 
E_n  & = &  \{(x,y,b)\in X \times Y\times \mathbb{R},\; f_n(x,y)\leq b\} \cap (X \times Y \times (-\infty,n]) \\
 X \times Y \times (-\infty,n] & \in &  \mathcal{B}(X\times Y\times \mathbb{R}) = \Delta_1^1(X \times Y \times \mathbb{R}) \subset \Delta_p^1(X\times Y\times \mathbb{R}),
\end{eqnarray*} see  (iii) and (iv) in Proposition \ref{base_hierarchy_pj}. 
 Then,  $(x,y,b) \mapsto f_n(x,y)$ is $\Delta_p^1(X \times Y \times \mathbb{R})$-measurable (see Proposition \ref{lemma_real_proj_pj} and 
 (v) in Proposition \ref{base_hierarchy_pj}). As 
  $(x,y,b) \mapsto b$ is $\mathcal{B}(X \times Y \times \mathbb{R})$-measurable, Proposition \ref{lemma_real_proj_pj} shows that $(x,y,b) \mapsto f_n(x,y)-b$ is $\Delta_p^1(X\times Y\times \mathbb{R})$-measurable and so, that $E_n \in \Delta_p^1(X\times Y \times \mathbb{R}),$ which concludes the proof of (i) when $f\geq 0$. 
In the general case, if $f : X \times Y \mapsto \mathbb{R} \cup \{-\infty, +\infty\}$, \eqref{base_integral_pj} implies that 
$$ \int_{-} f(x,y) q[dy|x] = \int f^+(x,y) q[dy|x] - \int f^-(x,y) q[dy|x],$$
using the convention $+\infty - \infty = -\infty + \infty = -\infty$. Proposition \ref{lemma_real_proj_pj} shows that $f^+ = \max(f,0)$ and $f^-=\max(-f,0)$ are $\Delta_{p}^1(X\times Y)$-measurable. So, $x \mapsto \int f^{\pm}(x,y) q[dy|x]$ are $\Delta_{p+r+2}^1(X)$-measurable and finally, $\lambda$ is $\Delta_{p+r+2}^1(X)$-measurable (see Proposition \ref{lemma_real_proj_pj} again). This shows (i).\\ 
(ii) If $f$ is a projective function and $q$ is a projectively measurable stochastic kernel, then there exist some $p\geq 1$ and $r\geq 1$ such that $f$ is $\Delta_p^1(X)$-measurable and $x\mapsto q[\cdot|x]$ is $\Delta_r^1(X)$-measurable and (i) shows (ii). \Halmos \\\endproof

We finish by showing that the multiple-priors conditional support of a projective function is a projective set-valued mapping. We first give the definition of measurability for set-valued mapping (see also \cite[Definition 14.1, p.643]{ref5_pj}). 
\begin{definition}\label{def_projective_mes_mapping}
Let $X$ be a Polish space and $d\geq 1$.  A set-valued mapping $F: X \twoheadrightarrow \mathbb{R}^d $ is $\Delta_n^1(X) $-measurable for some $n\geq 1$  if for every open sets $O \subseteq \mathbb{R}^d $ : 
$$F_{-1}(O) := \{ x \in X : F(x) \cap O \neq \emptyset \}  \in \Delta_n^1(X).$$  
The mapping $F $ is projective (or projectively measurable), if there exists $n \geq 1 $ such that $F $ is $\Delta_n^1(X) $-measurable.
\end{definition}
The next proposition  generalizes \citep[Lemma 4.3]{refbn1_pj} and  \citep[Lemma 2.6]{ref4_pj}. 
\begin{proposition}
Assume the (PD) axiom. Let $X$ and $Y$ be Polish spaces. 
Let  $\mathcal{Q}: X \twoheadrightarrow \mathfrak P(Y)$ be a nonempty set-valued mapping and  $f : X \times Y \to \mathbb{R}^d$ be a function. 
The multiple-priors conditional support set-valued mapping $D : X \twoheadrightarrow \mathbb{R}^d$ 
is defined by
\begin{eqnarray}
D(x)&:=&\bigcap \{A\subset \mathbb{R}^d,\; \mbox{closed},\; p[f(x,\cdot)\in A]=1,\;\forall p \in \mathcal{Q}(x)\} \label{support_def_pj}.
\end{eqnarray}
Let $r\geq 1$. Assume that  $\textup{Graph}(\mathcal{Q}):=\{(x,p)\in X \times \mathfrak P(Y), p\in\mathcal{Q}(x)\}\in \Delta_r^1(X \times \mathfrak P(Y))$ and $f$ is $\Delta_r^1(X\times Y)$-measurable. 
Then,  $D $ is a $\Delta_{r+1}^1(X)$-measurable set-valued mapping and $\textup{Graph}(\textup{Aff}(D )) \in \Delta_{r+1}^1( X) \otimes \mathcal{B}(\mathbb{R}^d)$. \\
Now, assume that  $\textup{Graph}(\mathcal{Q})$ is a projective set and $f$ is a projective function. Then, $D $  is a projective set-valued mapping and $\textup{Graph}(\textup{Aff}(D )) \in \textbf{P}( X \times \mathbb{R}^d)$. 
\label{aff_is_proj_pj}
\end{proposition}
\proof{Proof.}
For all $x\in X$, $h\in D (x)$ if and only if for all $n\geq 1$,  there exists $p\in\mathcal{Q} (x)$ such that $p[f(x,\cdot)\in \mathbb{B}(h,1/n)]>0$ where $\mathbb{B}(h,1/n)$ is the open ball of $\mathbb{R}^d$ centered at $h$ with radius $1/n$. So, for an open set $O$ of $\mathbb{R}^d$: 
\begin{eqnarray*}
D _{-1}(O)  := \{x\in  X,\; D (x) \cap O \neq \emptyset\}=\big\{x\in X,\; \exists h \in O,\; \forall n\geq 1,\; \exists p \in\mathcal{Q} (x),\; p\big[f(x,\cdot) \in \mathbb{B}(h,1/n)\big]) > 0\big\}.
\end{eqnarray*}
Let
\begin{eqnarray*}
F^n &:=& \big\{(x,h, y)\in X\times \mathbb{R}^d\times Y,\; f(x,y) \in \mathbb{B}(h,1/n)\big\}\\
F^n_{x,h}&:=& \big\{y\in Y,\; (x,h,y)\in F^n\big\} \;\; \mbox{for all } (x,h)\in X\times\mathbb{R}^d\\
E^n &:=& \big\{(x,h,p)\in X\times \mathbb{R}^d \times\mathfrak{P}(Y),\; p[F^n_{x,h}] > 0\big\}.
\end{eqnarray*}
Then, 
\begin{eqnarray}
D _{-1}(O) &=& \big\{x\in X,\; \exists h \in O,\; \forall n\geq 1,\; \exists p\in\mathcal{Q} (x),\; p[F^n_{x,h}]>0 \big\}\nonumber\\
&=& \textup{proj}_{ X}\big(( X\times O) \cap \bigcap_{n\geq 1} \{(x,h)\in X\times \mathbb{R}^d,\; \exists p\in \mathcal{Q} (x),\; p[F^n_{x,h}]>0 \}\big)\nonumber \\
&=& \textup{proj}_{ X}\Big(( X\times O) \cap \bigcap_{n\geq 1} \textup{proj}_{ X\times\mathbb{R}^d} \big(\sigma^{-1}(\textup{Graph}(\mathcal{Q} )\times \mathbb{R}^d)\cap E^n \big)\Big),\label{temp_AFF_0_pj}
\end{eqnarray}
where $\sigma(x,h,p)=(x,p,h)$ for all $(x,h,p)\in  X\times \mathbb{R}^d \times \mathfrak{P}(Y)$. \\
\noindent \textit{Assume that $\textup{Graph}(\mathcal{Q})\in \Delta_r^1(X \times \mathfrak P(Y))$ and $f$ is $\Delta_r^1(X\times Y)$-measurable. }\\
Assume for a moment that $E^n\in \Sigma_r^1( X\times \mathbb{R}^d\times \mathfrak{P}(Y))$. Recalling Proposition \ref{base_hierarchy_pj} (iii), (iv), (v),  and (i) as $\sigma$ is Borel, we have that 
\begin{eqnarray*} X \times O&  \in & \mathcal{B}( X\times \mathbb{R}^d) = \Delta_1^1( X\times \mathbb{R}^d) \subset \Sigma_1^1( X\times \mathbb{R}^d) \subset\Sigma_r^1( X\times \mathbb{R}^d)\\
\textup{Graph}(\mathcal{Q} ) \times \mathbb{R}^d & \in &  \Sigma_r^1\big( X\times \mathfrak{P}(Y)\big) \times \Sigma_r^1\big(\mathbb{R}^d) \subset \Sigma_r^1( X\times \mathfrak{P}(Y)\times \mathbb{R}^d\big)\\
\sigma^{-1}(\textup{Graph}(\mathcal{Q} )\times \mathbb{R}^d) & \in  & \Sigma_r^1( X\times \mathbb{R}^d\times \mathfrak{P}(Y)).
\end{eqnarray*}
So, recalling \eqref{temp_AFF_0_pj}, and using Proposition \ref{base_hierarchy_pj} (i) (stability of $\Sigma_r^1$ by projections and countable intersections), and then  \eqref{eq(v)_pj}, we find that  
$$D _{-1}(O) \in \Sigma_r^1( X) \subset \Delta_{r+1}^1( X).$$
Thus, the set-valued mapping $D $ is $\Delta_{r+1}^1( X)$-measurable. Applying successively \citep[Exercise 14.12, p652]{ref5_pj} and \citep[Theorem 14.8, p648]{ref5_pj}, we get that $\textup{Graph}(\textup{Aff}(D )) \in \Delta_{r+1}^1( X) \otimes \mathcal{B}(\mathbb{R}^d)$. \\
So, it remains to show that $E^n \in \Sigma_r^1( X\times \mathbb{R}^d \times \mathfrak{P}(Y))$. For all $q\geq 1$, let
$$E_q^n := \big\{(x,h,p) \in  X \times \mathbb{R}^d \times \mathfrak{P}( Y),\; p[F_{x,h}^n]\geq 1/q\big\}.$$
Then, $E^n = \bigcup_{q\geq 1} E_q^n$. If we show that for all $q\geq 1$, $E^n_q\in \Sigma_r^1( X \times \mathbb{R}^d \times \mathfrak{P}(Y))$, Proposition \ref{base_hierarchy_pj} (i) proves the claim. First, remark that the function $(x,h,y) \mapsto f(x,y)$ is $\Delta_r^1( X\times \mathbb{R}^d \times Y)$-measurable, see Lemma \ref{etendre_fct_mes}. So, 
$$F^n \in \Delta_r^1( X\times \mathbb{R}^d\times Y)\subset \Sigma_r^1( X\times \mathbb{R}^d\times Y).$$
Now, we have using Fubini's theorem that $p[F_{x,h}^n] = (\delta_{x,h} \otimes p)[F^n]$ for all $(x,h,p) \in  X\times \mathbb{R}^d \times\mathfrak{P}(Y)$. Let $q\geq 1$. Then, 
\begin{eqnarray*}
E^n_q &=& \kappa^{-1}\big(\{\nu\in\mathfrak{P}( X \times \mathbb{R}^d \times Y),\; \nu[F^n] \geq 1/q\}\big),
\end{eqnarray*}
where $\kappa :  X \times \mathbb{R}^d \times \mathfrak{P}(Y) \to \mathfrak{P}( X\times \mathbb{R}^d \times Y)$ is such that $\kappa(x,h,p)=\delta_{x,h} \otimes p.$ 
Then, $\kappa$ is continuous (see \citep [Lemma 7.12, p144]{ref1_pj} and \citep [Corollary 7.21.1, p130]{ref1_pj}). 
Recalling that $F^n\in\Sigma_r^1( X\times \mathbb{R}^d\times Y)$ and using Proposition \ref{prop_crutial_hierarchy_pj}, we obtain that 
$$\{\nu\in\mathfrak{P}( X \times \mathbb{R}^d \times Y),\; \nu(F^n) \geq 1/q\}\in \Sigma_r^1(\mathfrak{P}( X \times \mathbb{R}^d \times Y)).$$ 
Thus, (i) in Proposition \ref{base_hierarchy_pj} implies that $E^n_q \in \Sigma_r^1( X\times \mathbb{R}^d\times \mathfrak{P}(Y))$.\\
\noindent \textit{Assume that  $\textup{Graph}(\mathcal{Q})$ is a projective set and $f$ is a projective function.  }\\
There exist some $q, r\geq 1$ such that $\textup{Graph}(\mathcal{Q} )\in\Delta_q^1( X\times \mathfrak{P}(Y))$ and $f$ is $\Delta_r^1(X\times Y)$-measurable.  We may assume that $q\leq r$. Then, $f$ is $\Delta_r^1(X \times Y)$-measurable, see Proposition \ref{base_hierarchy_pj} (iv). So,  $D $ is a  $\Delta_{r+1}^1(X)$-measurable set-valued mapping and $\textup{Graph}(\textup{Aff}(D )) \in \Delta_{r+1}^1( X) \otimes \mathcal{B}(\mathbb{R}^d)$. It follows that $D$ is a projective set-valued mapping. 
Now, using (iii) to (v) in Proposition \ref{base_hierarchy_pj}, 
$$\textup{Graph}(\textup{Aff}(D )) \in\Delta_{r+1}^1( X) \otimes \mathcal{B}(\mathbb{R}^d)\subset \Delta_{r+1}^1( X)\times \Delta_{r+1}^1(\mathbb{R}^d) \subset \Delta_{r+1}^1( X\times \mathbb{R}^d) \subset \textbf{P}( X\times \mathbb{R}^d),$$ 
and the proof is completed. 
\Halmos \\ \endproof

\section{Applications to mathematical finance and economic}
\label{secappli}
Traditional approaches in mathematical finance and economic assume a single probability measure $P$ to model the law of the asset prices. 
Knightian uncertainty (name after the work of Knight, see \cite{refkni_pj}) allows that agents may hold diverse and even conflicting beliefs about future states of the world and consider a set $\mathcal{Q}$ of probability measures (or priors)  instead of a single one. 
The earliest literature assumed that $\mathcal{Q}$ is dominated by one probability measure. We refer to \cite{refdoml_pj} for a comprehensive survey of the dominated case. Unfortunately, this setting excludes volatility uncertainty and is easily violated in discrete time (see \cite{ref4_pj}); this is why the latest literature has focused on the non-dominated case.   
Bouchard and Nutz have developed a setting for the so-called quasi-sure non-dominated uncertainty in  \citep{refbn1_pj}. The usual notion of $P$-almost-sure is replaced by that of $\mathcal{Q}$-quasi-sure, where loosely speaking, a property holds $\mathcal{Q}$-quasi-surely if it holds $P$-almost-surely, for every $P\in \mathcal{Q}$. 
Set-valued mappings of ``local" priors (investor's beliefs between two successive moments) are first given.  The cornerstone assumption of \cite{refbn1_pj} is that the graphs of these set-valued mappings are analytic sets. 
Thanks to this assumption and the measurable selection theorem of Jankov-von Neumann, it is possible to obtain local beliefs that are analytically and, thus, universally measurable as a function of the path.  The inter-temporal set of beliefs can then be constructed from these kernels as product measures. Measurable selection is also necessary to do the way back, for example, to go from inter-temporal quasi-sure inequalities to local quasi-sure ones. 
For that, Bouchard and Nutz rely on the uniformisation of nuclei of Souslin schemes for the product of the universal sigma-algebra and the Borel one, as discussed by Leese in \cite{Leese78}.  So, already in \citep{refbn1_pj}, one needs to go outside the class of analytic sets, which is the class of all nuclei of Souslin schemes for the Borel sigma-algebra.  
Moreover, Bouchard and Nutz use upper-semianalytic functions as a measurability standard. As already mentioned, the class of upper-semianalytic functions is not closed under composition. However, the composition of an upper-semianalytic function and a Borel measurable one remains upper-semianalytic. This is why they assume the asset prices to be Borel measurable.  
 Furthermore, as the complement of an analytic set may not be analytic, this generates further measurability issues in  \citep{refbn1_pj} (for example, they have to restrict upper-semianalytic functions to coanalytic sets). 
Summing up in the classical framework of Bouchard and Nutz, the price processes are assumed to be Borel measurable, the graphs of the local random beliefs to be analytic sets, while trading strategies and stochastic kernels are only obtained to be universally measurable by measurable selection.  

 In the companion paper \cite{cf2024}, we have introduced the projective setting, which unifies these measurability conditions. The prices $S_t$ are assumed to be projective functions, the graphs of the local beliefs set-valued mappings to be projective sets, while trading strategies and stochastic kernels are obtained to be projectively measurable by measurable selection, see Proposition \ref{univ_sel_proj_pj}.  
So, \cite{cf2024} proposes a homogeneous setup for Knightian uncertainty, which is more general than that of Bouchard and Nutz. Indeed, analytic sets are projective sets, while Borel measurable functions are projective ones without postulating the (PD) axiom. 
In this projective setup, \cite{cf2024} proves, under the (PD) axiom,  the existence of an optimal investment strategy for the following maxmin problem, based on the representation of preferences of \cite{refgs_pj}:
\begin{eqnarray*}
u(x)=\sup_{\phi}\inf_{P\in \mathcal{Q}}\mathbb{E}_P U(\cdot,V_T^{x,\phi}(\cdot)).
\label{RUMP}
\end{eqnarray*} 
where $U$ is a projective function, non-decreasing, upper semicontinuous but not concave in $y\in \mathbb{R}$, $\phi$ is the control (the stochastic process called portfolio or strategy in the financial context) and $V_T^{x,\phi}=x+\sum_{t=1}^T <\phi_t,S_t>$ is the wealth obtained with the control $\phi$ at the terminal time $T$, starting from an initial wealth $x$. 
In  \cite{cf2024}, we also provide a counter-example, where we show that ZFC is not sufficient for solving a multiperiod optimisation problem, when $U$ is less than lower-semianalytic. Indeed, we show that in the constructible universe $L$, one can find a Borel price process, priors which graphs are Borel sets, where the utility function is bounded from above, increasing, concave and continuous as of function of $y$, and only analytically measurable (i.e., measurable with respect to the sigma-algebra generated by the analytic sets), as a function of the path, and where the unique optimal strategy is projective, but not Lebesgue measurable. 
However, the constructible universe $L$ is inconsistent with the (PD) axiom, since the axiom of Constructibility holds in $L$. Having a set of admissible strategies, which includes non-Lebesgue measurable functions, leads to an ill-posed problem for manipulating integrals of such functions, which are not well-defined without the (PD) axiom. This justifies using models with ZFC+(PD).  

Later, in the projective setting, \cite{bci25} provides characterisations of the quasi-sure No-Arbitrage condition. This condition states that if the terminal wealth starting from 0, is non-negative $\mathcal{Q}$-quasi-surely, then it equals $0$ $\mathcal{Q}$-quasi-surely, i.e., 
$$
\quad V_T^{0, \phi} \geq 0 \; \mathcal{Q} \mbox{-q.s.} \text{ for some }  \phi  \implies V_T^{0, \phi} = 0 \; \mathcal{Q}\mbox{-q.s.}
$$

Both \cite{cf2024} and \cite{bci25} use the results that we prove here. Not only are the results on projective sets and functions used, but also the results on $\Delta_n^1(X)$-measurable functions and on the other classes of sets of the hierarchy. Indeed, a lot of theorems of analysis require working with a sigma-algebra. This is the case, for example,  in \cite[Theorem 18.5]{refAB_pj}, which shows that the distance between a point and a set-valued mapping is a Caratheodory function and is crucially used in \cite{bci25}. We have also used this for applying  \citep[Theorem 14.8, p648]{ref5_pj}, in the proof of Proposition \ref{aff_is_proj_pj}.  
Moreover, it is often necessary in \cite{cf2024,bci25} to use countable supremum. The class of projective functions is not closed under countable supremum, but the class of $\Delta^1_n$-measurable functions is. For more technical details,  we refer, for example, to \cite[Remark 6]{bci25}.

So, the present paper provides the mathematical basis for developing more models that can capture the complexities of uncertainty without the traps coming from measurable selections. The interplay between projective sets, projective functions, and the (PD) axiom provides a powerful foundation for understanding dynamic decision-making in ambiguous environments, reinforcing the theoretical underpinnings of multiple-priors financial models.

\section{Appendix}
\label{appendix_pj}

\label{hiera_section_proof}
{\bf Proof of Proposition \ref{base_hierarchy_pj}}\\
First, the equality $\Delta_1^1(X) = \mathcal{B}(X)$  in (iii) is proved for example in \citep[Theorem 14.11, p88]{refproj_pj}.  

\noindent  \textit{Proof of (iv).}\\
We show by induction that for all $n\geq 1$ and for all Polish spaces $X$, $\Sigma_{n}^1(X) \subset \Sigma_{n+1}^1(X)$. \\
We start with $n=1$. 
Let $X$ be a Polish space and $A\in \Sigma_{1}^1(X)$ be an analytic set of $X$. 
There exists some $C\in \mathcal{B}(X\times \mathcal{N})$ such that $A = \textup{proj}_X(C)$, see \eqref{level1}. As $X\times \mathcal{N}$ is a Polish space, $\mathcal{B}(X \times \mathcal{N}) = \Delta_1^1(X \times \mathcal{N})\subset \Pi_1^1(X \times \mathcal{N})$, see \eqref{Delta_pj}, and $C\in \Pi_1^1(X \times \mathcal{N})$. So, $A\in \Sigma_{2}^1(X)$. This shows the initialisation step. \\
Assume now that the induction hypothesis holds for $n\geq 1$. Let $X$ be a Polish space and $A \in \Sigma_{n+1}^1(X)$. Then, there exists some $C\in \Pi_n^1(X\times \mathcal{N})$ such that $A = \textup{proj}_X(C)$. Using the induction hypothesis for $X\times \mathcal{N}$, which is a Polish space, we get that $\Sigma_n^1(X\times \mathcal{N})\subset \Sigma_{n+1}^1(X\times \mathcal{N})$ and \eqref{Pi_pj} implies also the inclusion $\Pi_{n}^1(X\times \mathcal{N}) \subset \Pi_{n+1}^1(X\times \mathcal{N})$. 
Thus, $C \in \Pi_{n+1}^1(X\times \mathcal{N})$ and $A \in \Sigma_{n+2}^1(X)$. This concludes the induction and the sequence $(\Sigma_n^1(X))_{n\geq 1}$ is nondecreasing.  Then, we trivially have, using  \eqref{Pi_pj} and \eqref{Delta_pj}, that the sequences $(\Pi_n^1(X))_{n\geq 1}$ and $(\Delta_n^1(X))_{n\geq 1}$ are nondecreasing.

\noindent  \textit{Proof of (i), (ii) and (iii).}\\
We show by induction that for all $n\geq 1$ and for all Polish spaces $X$, (i), (ii) and (iii) hold.  \\
We start with $n=1$. 
Let $X$ be a Polish space. Assertion (i) for $n=1$ follows from \citep[Corollary 7.35.2, p160]{ref1_pj} and \citep[Proposition 7.40, p165]{ref1_pj}. Moreover, (ii) and (iii)  follow easily from (i), \eqref{Pi_pj}, \eqref{Delta_pj} and from  $\mathcal{B}(X)=\Delta_1^1(X)$. This shows the initialisation step. \\
Assume now that for all Polish spaces $X$, (i), (ii) and (iii) hold for $n\geq 1$. 
Assume for a moment that (i) holds for $n+1$. Then, the closeness by countable unions and intersections and by preimage in (ii) and (iii) is a direct consequence of (i), as in the initialisation step. To prove that $\Delta_{n+1}^1(X)$ is a sigma-algebra, it remains to prove the closeness by complement. Let $A\in\Delta_{n+1}^1(X) = \Sigma_{n+1}^1(X) \cap \Pi_{n+1}^1(X)$. Using \eqref{Pi_pj}, $A\in\Sigma_{n+1}^1(X)$ implies that $X \setminus A \in \Pi_{n+1}^1(X)$ and $A\in\Pi_{n+1}^1(X)$ implies that $X \setminus A \in \Sigma_{n+1}^1(X)$. Thus, $X \setminus A \in \Delta_{n+1}^1(X)$. It remains to prove (i).

\textit{Proof of (i) for $n+1$ : Countable unions.}\\
Let $X$ be a Polish space. Let $(A_i)_{i\geq 0} \subset \Sigma_{n+1}^1(X)$ and $(C_i)_{i\geq 0} \subset \Pi_{n}^1(X\times \mathcal{N})$ be such that $A_i=\textup{proj}_X(C_i)$ for all $i\geq 0$. We have that 
\begin{eqnarray*}
x\in \bigcup_{i \geq 0} A_i \iff \exists i \geq 0,\; x\in A_i &\iff &  \exists i \geq 0,\; \exists u \in \mathcal{N},\;  (x,u)\in C_i \\ 
&\iff & \exists u \in \mathcal{N},\; (x,u)\in \bigcup_{i\geq 0}C_i \iff x\in \textup{proj}_X\bigl(\bigcup_{i\geq 0}C_i\bigr).
\end{eqnarray*}
 Now, (ii) for $n$ and $X\times \mathcal{N}$ shows that $\cup_{i\geq 0}C_i \in \Pi_n^1(X\times \mathcal{N})$. Thus, using \eqref{Sigma_pj}
 $$\bigcup_{i\geq 0} A_i = \textup{proj}_X\bigl(\bigcup_{i\geq 0}C_i\bigr)\in \Sigma_{n+1}^1(X).$$

\textit{Proof of (i) for $n+1$ : Borel preimages.}\\
Let $X$ and $Y$ be Polish spaces. Let $f : X \to Y$ be Borel measurable. Let $B \in \Sigma_{n+1}^1(Y)$ and $C\in \Pi_n^1(Y\times \mathcal{N})$ such that $B= \textup{proj}_Y(C)$. Then, 
\begin{eqnarray*}
f^{-1}(B) 
&=&\{x\in X,\; \exists u\in\mathcal{N},\; (f(x),u) \in C \} = 
\{x\in X,\; \exists u\in\mathcal{N},\; (x,u) \in \Psi^{-1}(C)\}=\textup{proj}_X(\Psi^{-1}(C)), 
\end{eqnarray*}
where $ X\times\mathcal{N} \ni (x,u) \mapsto\Psi(x,u):=(f(x),u)\in Y\times \mathcal{N}$ is Borel measurable.  Using (ii) for $n$ and $X\times \mathcal{N}$, we have that $\Psi^{-1}(C) \in \Pi_n^1(X\times \mathcal{N}).$ So, \eqref{Sigma_pj} implies that $f^{-1}(B) \in \Sigma_{n+1}^1(X)$.

\textit{Proof of (i) for $n+1$ : Borel images.}\\
Let $X$ and $Y$ be Polish spaces. Let $f : X \to Y$ be Borel measurable. Let $A \in \Sigma_{n+1}^1(X)$ and $C\in \Pi_n^1(X\times \mathcal{N})$ such that $A= \textup{proj}_X(C)$. As $X\times\mathcal{N}$ is a non-empty Polish space, \citep[Theorem 7.9, p38]{refproj_pj} shows that there exists a continuous surjection $\Phi$ from $\mathcal{N}$ to $X\times \mathcal{N}$. Then,
\begin{eqnarray*}
f(A) &=& \{y\in Y,\; \exists x \in A,\; y=f(x)\}= \{y\in Y,\; \exists (x,u) \in X\times\mathcal{N},\; y=f(x) \;\; \mbox{and} \;\; (x,u)\in C\}\\ 
&=& \{y\in Y,\; \exists w \in \mathcal{N},\; y=g(w) \;\; \mbox{and} \;\; \Phi(w)\in C\}= \textup{proj}_Y\bigl(\xi^{-1}(\textup{Graph}(g)) \cap (Y \times \Phi^{-1}(C))\bigr),
\end{eqnarray*}
where for all $w\in \mathcal{N}$, $g(w) := f(\textup{proj}_X(\Phi(w)))$ and for all $(y,w)\in Y\times \mathcal{N}$, $\xi(y,w):=(w,y)$. \\
Assume for a moment that both $Y\times \Phi^{-1}(C)$ and $\xi^{-1}(\textup{Graph}(g))$ belong to $\Pi_n^1(Y\times \mathcal{N}).$ Then (ii) for $n$ and $X\times \mathcal{N}$ shows that their intersection remains in $\Pi_n^1(Y\times \mathcal{N}).$ So, using \eqref{Sigma_pj},  $f(A)\in \Sigma_{n+1}^1(Y)$. \\
Note first that $\Phi^{-1}(C) \in \Pi_n^1(\mathcal{N})$ using (ii) for $n$ as $C \in \Pi_n^1(X \times \mathcal{N})$ and $\Phi$ is continuous and thus Borel measurable. Now, as $Y \times \mathcal{N}\ni(a,b)\mapsto b$ is Borel measurable, (ii) for $n$ implies that $Y\times \Phi^{-1}(C) \in \Pi_n^1(Y\times \mathcal{N})$. \\
Finally, as $g$ is Borel measurable, \citep[Corollary 7.14.1, p121]{ref1_pj} shows that 
$$\textup{Graph}(g) \in \mathcal{B}(\mathcal{N}\times Y) = \Delta_1^1(\mathcal{N}\times Y)\subset \Pi_1^1(\mathcal{N}\times Y).$$
But (iv) shows that $\Pi_1^1(\mathcal{N}\times Y) \subset \Pi_n^1(\mathcal{N}\times Y)$ and $\textup{Graph}(g) \in \Pi_n^1(\mathcal{N}\times Y)$. Using again (ii) for $n$ as $\xi$ is Borel measurable, we deduce that $\xi^{-1}(\textup{Graph}(g))\in \Pi_n^1(Y\times \mathcal{N})$.

\textit{Proof of (i) for $n+1$ : Countable intersections.}\\
Let $X$ be a Polish space. Let $(A_i)_{i\geq 0} \subset \Sigma_{n+1}^1(X)$ and $(C_i)_{i\geq 0} \subset \Pi_{n}^1(X\times \mathcal{N})$ be such that $A_i=\textup{proj}_X(C_i)$ for all $i\geq 0$. 
Then,
\begin{eqnarray*}
\bigcap_{i\geq 0} A_i &=& \{x\in X,\; \forall i \geq 0,\; \exists u_i\in\mathcal{N},\; (x,u_i)\in C_i\}\\
&=&\{x\in X,\; \exists u=(u_i)_{i\geq 0}\in \mathcal{N}^\mathbb{N},\; \forall i\geq 0,\; (x,u_i)\in C_i\}\\
&=&\{x\in X,\; \exists z\in\mathcal{N},\; \forall i\geq 0,\; (x,\Phi(z)_i)\in C_i\}= \textup{proj}_X\bigl(\bigcap_{i\geq 0}D_i\bigr),
\end{eqnarray*}
where $\Phi$ is a continuous surjection from $\mathcal{N}$ to the Polish space $\mathcal{N}^\mathbb{N}$ (see \citep[Proposition 7.4, p108]{ref1_pj}) given by \citep[Theorem 7.9, p38]{refproj_pj} and for all $i\geq 0$, $D_i:= \{(x,z)\in X \times \mathcal{N},\; (x, \Phi(z)_i) \in C_i\}.$
As for all $i\geq 0$, $(x,z) \mapsto (x, \Phi(z)_i)$ is Borel measurable and $C_i\in\Pi_n^1(X\times\mathcal{N)}$, (ii) for $n$ shows that $D_i\in \Pi_{n}^1(X\times \mathcal{N})$. So, again (ii) for $n$ and $X\times \mathcal{N}$ shows that $\bigcap_{i\geq 0}D_i \in \Pi_{n}^1(X\times \mathcal{N})$. Thus, $\bigcap_{i\geq 0} A_i \in \Sigma_{n+1}^1(X)$, see \eqref{Sigma_pj}. This concludes the induction for (i).

\noindent \textit{Proof of (v).}\\
Let $n\geq 1$ and $X$ be a Polish space. We show the first inclusion in \eqref{eq(v)_pj}.\\ 
For that, we prove that $\Pi_n^1(X) \subset \Sigma_{n+1}^1(X)$. Let $A\in\Pi_n^1(X).$ As $X\times \mathcal{N} \ni (a,b) \mapsto a$ is Borel measurable, (ii) for $n$ and $X\times \mathcal{N}$ implies that $A\times \mathcal{N}\in \Pi_n^1(X\times \mathcal{N})$. Thus, $A = \textup{proj}_{X}(A\times \mathcal{N})\in \Sigma_{n+1}^1(X)$, see \eqref{Sigma_pj}. \\
Now, let $B\in \Sigma_n^1(X)$. Then, using \eqref{Pi_pj}, $X\setminus B \in \Pi_n^1(X) \subset \Sigma_{n+1}^1(X)$ and so, $B\in \Pi_{n+1}^1(X)$.   
This proves that  $\Sigma_n^1(X) \subset \Pi_{n+1}^1(X)$. Recalling (iv), we also have that $\Sigma_{n}^1(X)\subset \Sigma_{n+1}^1(X)$, and thus 
$$\Sigma_{n}^1(X) \subset \Pi_{n+1}^1(X) \cap \Sigma_{n+1}^1(X)=\Delta_{n+1}^1(X),$$ 
using \eqref{Delta_pj}. 
As $\Delta_{n+1}^1(X)$ is a sigma-algebra, \eqref{Pi_pj} proves that $\Pi_{n}^1(X)\subset \Delta_{n+1}^1(X)$. This shows that $\Sigma_{n}^1(X) \cup \Pi_{n}^1(X) \subset \Delta_{n+1}^1(X)$. 
Using this last inclusion and \eqref{Delta_pj}, we have the following inclusions for all $n\geq 1$:
$$\Pi_n^1(X) \subset \Delta_{n+1}^1(X) \subset \Sigma_{n+1}^1(X) \subset \Delta_{n+2}^1(X) \subset \Pi_{n+2}^1(X).$$ 
Thus, taking the union over $n\geq 1$ and using (iv), we get that the unions of $(\Sigma_n^1(X))_{n\geq 1}$, $(\Pi_n^1(X))_{n\geq 1}$ and $(\Delta_n^1(X))_{n\geq 1}$ are equal, which achieves the proof of \eqref{eq(v)_pj}. \\
We now prove that $\Sigma_n^1(X) \times \Sigma_n^1(Y) \subset \Sigma_n^1(X \times Y)$. Let $A \in\Sigma_n^1(X)$ and $B\in\Sigma_n^1(Y)$. As the mappings $X \times Y\ni(a,b) \mapsto a $ and $X \times Y \ni (a,b)\mapsto b $ are Borel measurable, (i) shows that  $A\times Y$ and $X \times B$ belong to $\Sigma_n^1(X\times Y)$ and also that $$A \times B =  (A\times Y) \cap  (X \times B)\in \Sigma_n^1(X\times Y).$$ 
The same method, using (ii) and (iii) instead of (i), shows that $\Pi_n^1(X) \times \Pi_n^1(Y) \subset \Pi_n^1(X\times Y)$ and $\Delta_n^1(X) \times \Delta_n^1(Y) \subset \Delta_n^1(X\times Y)$.

\noindent \textit{Proof of (vi).}\\
We first show that  $\textbf{P}(X)$ is closed under complement. Let $A \in \textbf{P}(X)$. Then, there exists some $n\geq 1$ such that $A \in\Delta_n^1(X)$. As $\Delta_n^1(X)$ is a sigma-algebra, $X \setminus A \in \Delta_n^1(X) \subset \textbf{P}(X)$. We now prove that $\textbf{P}(X)$ is closed under finite unions and intersections. Let $A \in \textbf{P}(X)$ and $B\in\textbf{P}(X)$, there exist $n, p \geq 1$ such that $A\in\Delta_n^1(X)$ and $B\in\Delta_p^1(X)$. We may assume that $n\leq p$. Then, (iv) and (iii) show that $A,$ $A\cup B$ and $A \cap B$ belong to $\Delta_p^1(X) \subset \textbf{P}(X)$.\\ Now, let $f : X \to Y$ be Borel measurable. Let $B \in \textbf{P}(Y)$. Then, there exists some $q\geq 1$ such that $B \in \Delta^1_q(Y)$ and it follows from (iii) that $f^{-1}(B)\in\Delta_q^1(X) \subset \textbf{P}(X)$. Let $A \in \textbf{P}(X)$. Then, there exists some $q\geq 1$ such that $A \in \Delta^1_q(X) \subset \Sigma_{q}^1(X)$ (see \eqref{Delta_pj}) and it follows from (i) and  \eqref{eq(v)_pj} that 
$$f(A) \in \Sigma_q^1(Y) \subset \Delta_{q+1}^1(Y) \subset \textbf{P}(Y).$$
Finally, we now prove \eqref{eq(vi)_pj}. 
As $\Delta_1^1(X) = \mathcal{B}(X)$, $\mathcal{B}(X) \subset \textbf{P}(X)$. Now, \eqref{eq(v)_pj} shows that 
$$\Sigma_1^1(X)  \cup \Pi_1^1(X) \subset \Delta_2^1(X) \subset \textbf{P}(X).$$ 
Let $A \in \textbf{P}(X)$ and $B\in \textbf{P}(Y)$. 
There exist $n,p \geq 1$ such that $A \in \Delta_n^1(X)$ and $B\in\Delta_p^1(Y)$. Again, we may assume that $n\leq p$ and we have that $A\in\Delta_p^1(X)$. Using (v), we get that $$A \times B \in \Delta_{p}^1(X\times Y) \subset \textbf{P}(X\times Y). \Halmos$$  \endproof

\noindent {\bf Proof of Proposition \ref{lemma_real_proj_pj}}\\
(i) The first part of the proof follows from the fact that $\Delta_p^1(X)$ is a sigma-algebra (see Proposition \ref{base_hierarchy_pj} (iii)) and from conventions \eqref{cvt_inf_pj} and \eqref{cvt_new} and is given for the sake of completeness. 
We first prove the sup and inf results. For all $c\in\mathbb{R},$ 
$$\{\min(f,g)< c\} = \{f< c\} \cup \{g< c\} \in \Delta_p^1(X) \mbox{ and } \{\max(f,g)< c\} = \{f < c\} \cap \{g<c\} \in \Delta_p^1(X).$$ 
So, $\min(f,g)$ and $\max(f,g)$ are $\Delta_p^1(X)$-measurable. Similarly, 
$$\bigl\{\sup_{n\geq 0} f_n \leq  c\bigr\} = \bigcap_{n\geq 0}\{f_n \leq c\} \in \Delta_p^1(X) \mbox{ and } \bigl\{\inf_{n\geq 0} f_n < c\bigr\} = \bigcup_{n\geq 0}\{f_n  < c\} \in \Delta_p^1(X)$$ and $\sup_{n\geq 0} f_n$ and $\inf_{n\geq 0} f_n$ are also $\Delta_p^1(X)$-measurable.  As $\{-f < c\} = f^{-1}((-c,+\infty]) \in \Delta_{p}^1(X)$, we have that $-f$ is $\Delta_p^1(X)$-measurable. 
We now show that $f+g$ is $\Delta_p^1(X)$-measurable. Let 
$$N := \{f=+\infty, g=-\infty\} \cup \{f=-\infty, g=+\infty\}.$$
Then, as $ \{f=+\infty\}=\cap_{n\geq 0}\{f \geq n\}$ and $ \{f=-\infty\}=\cap_{n\geq 0}\{f \leq -n\},$ $N\in\Delta_p^1(X)$.  Now, for all $c\in\mathbb{R}$, recalling convention \eqref{cvt_inf_pj}, 
$$\{x\in X,\; f(x)+g(x) < c\} = N \cup \Bigl((X\setminus N) \cap \bigcup_{r\in\mathbb{Q}}\{x\in X, f(x)<r\} \cap \{x\in X, g(x)<c-r\}\Bigr) \in \Delta_p^1(X).$$
Thus, $f+g$ is $\Delta_p^1(X)$-measurable. Now, we prove that $fg$ is $\Delta_p^1(X)$-measurable. We start with the case where $f\geq 0$ and $g\geq 0$. We have that $\{fg <c\} = \emptyset \in \Delta_p^1(X)$ for all $c\leq 0$. Due to convention \eqref{cvt_new}, we have that for all $c> 0$,  
$$\{x\in X,\; f(x)g(x) < c\} = \{x\in X,\; g(x)=0\} \cup \bigcup_{r\in\mathbb{Q}^{*}_+} \{x\in X, f(x)<c/r\} \cap \{x\in X, g(x)<r\} \in \Delta_p^1(X).$$
So, $fg$ is $\Delta_p^1(X)$-measurable when $f\geq 0$ and $g\geq 0$. The general case follows from the equality 
$$fg = (f^+ g^+ + f^- g^-) - (f^+ g^- + f^- g^+)$$  
and the fact that the sum (and the difference) of $\Delta_p^1(X)$-measurable functions are $\Delta_p^1(X)$-measurable. Assume now that $f> 0$. Then, $\{f^a< c\} = \emptyset \in \Delta_p^1(X)$ for all $c\leq 0$. If  $a>0$, $\{f^a< c\} = \{f< c^{1/a}\} \in\Delta_p^1(X)$ for all $c>0$. Now, if $a<0$, $\{f^a< c\} = \{f> c^{1/a}\} \in\Delta_p^1(X)$ for all $c>0$. So, $f^a$ is also $\Delta_p^1(X)$-measurable.\\ 
(ii)  The proof is the same as in Lemma \ref{couple_fct_mes} (ii). 
If $f$ and $g$ are projective functions, there exist $n, p \geq 1$ such that $f$ is $\Delta_n^1(X)$-measurable and $g$ is $\Delta_p^1(X)$-measurable. We may assume that $n\leq p$. Then, Proposition \ref{base_hierarchy_pj} (iv) shows that $f$ is $\Delta_p^1(X)$-measurable. So, (i) shows that $fg$, $f+g$, $-f$, $\min(f,g)$, $\max(f,g)$ and $f^a$ (if $f> 0$) are $\Delta_p^1(X)$-measurable and thus projectively measurable.
\Halmos \\ \endproof

\end{document}